\documentclass[11pt]{article}
\usepackage{epsfig,amsmath,amsthm, wasysym, xcolor}
\usepackage{algorithm}
\usepackage[noend]{algorithmic}
\usepackage{graphicx, tikz, array}

\usepackage{hyperref}

\newtheorem{df}{Definition}[section]
\newtheorem{definition}[df]{Definition}
\newtheorem{theorem}[df]{Theorem}
\newtheorem{lemma}[df]{Lemma}
\newtheorem{remark}[df]{Remark}

\newtheorem{corollary}[df]{Corollary}

\newtheorem{proposition}[df]{Proposition}
\usepackage{tcolorbox}

\title{Tyshkevich's Graph Decomposition and \\the Distinguishing Numbers of Unigraphs}
\author{Christine T. Cheng}
\date{Department of Electrical Engineering and Computer Science \\ University of Wisconsin-Milwaukee\\ ccheng@uwm.edu%
       }

\begin{document}
\maketitle

\begin{abstract}
A  $c$-labeling $\phi: V(G) \rightarrow \{1, 2, \hdots, c \}$ of graph $G$ is {\it distinguishing} if,  for every non-trivial automorphism $\pi$ of $G$,  there is some vertex $v$ so that $\phi(v) \neq \phi(\pi(v))$. The {\it distinguishing number of $G$}, $D(G)$, is the smallest $c$ such that $G$ has a distinguishing $c$-labeling.

We consider a compact version of Tyshkevich's graph decomposition theorem where trivial components are maximally combined to form a complete graph or a graph of isolated vertices.  Suppose the compact canonical decomposition of $G$ is  $G_{k} \circ G_{k-1} \circ  \cdots \circ G_1 \circ G_0$.  We prove that $\phi$ is a distinguishing labeling of $G$ if and only if $\phi$ is a distinguishing labeling of $G_i$ when restricted to $V(G_i)$ for $i = 0, \hdots, k$.  Thus, $D(G) = \max \{D(G_i), i = 0, \hdots, k \}$. We then present an algorithm that computes the distinguishing number of a unigraph in linear time.

\end{abstract}

\section{Introduction}

	A  $c$-labeling $\phi: V(G) \rightarrow \{1, 2, \hdots, c \}$ of graph $G$ is {\it distinguishing} if,  for every non-trivial automorphism $\pi$ of $G$,  there is some vertex $v$ so that $\phi(v) \neq \phi(\pi(v))$.  That is, $\phi$ breaks all the non-trivial symmetries of $G$ so that the only automorphism of the labeled graph $(G, \phi)$ is the identity map. The {\it distinguishing number of $G$}, $D(G)$, is the smallest $c$ such that $G$ has a distinguishing $c$-labeling.   For example, $D(K_n) = n$, $D(P_n) = 2$ for $n \ge 2$ and $D(K_{n, n}) = n+1$.   
	
	Babai used distinguishing labelings, calling them {\it asymmetric colorings}, to create asymmetric infinite trees with two prescribed degrees in his 1977 paper \cite{Ba77}.  But it was Albertson and Collins \cite{AlCo96} who introduced the notion of distinguishing labelings and numbers for {\it all} graphs in 1996.  Since then hundreds of papers in graph theory and group theory have been published on this topic (e.g. \cite{KlWoZh06}, \cite{ ImSmTuWa15}, \cite{CoTr22},  \cite{ShAhTaHa22}).

	Of interest to us are the algorithmic aspects of distinguishing labelings and numbers.  Russell and Sundaram \cite{RuSu98} proved that $DIST = \{(G,c)| \mbox{$G$ has a distinguishing $c$-labeling} \}$  is in the complexity class AM, the set of languages for which there are Arthur and Merlin games.   If restricted to certain graph families however, $DIST$ becomes tractable.  For example, efficient algorithms have been designed to compute the distinguishing numbers of  trees and forests \cite{ArDe04, Ch06}, planar graphs \cite{ArDe04, ArChDe08} and interval graphs \cite{Ch09}.  All of them are  based on efficient graph isomorphism algorithms  that were modified to count the number of inequivalent distinguishing labelings of a graph. This observation led us to {\it unigraphs} which arguably has one of the simplest graph isomorphism algorithms:  given a unigraph $G$ and another graph $H$,  $G \cong H$  if and only if $G$ and $H$ have identical degree sequences.  
We want to know if we can leverage this isomorphism algorithm to compute the distinguishing number of a unigraph.  As we shall show, the answer is a partial yes.  As is, the algorithm is too simple and does not reveal much information about the automorphisms of the unigraph. We need to rely on Tyshkevich's decomposition algorithm for {\it all} graphs.   Once we have a unigraph's decomposition, only then are we able to take advantage of the isomorphism algorithm to identify the components in the decomposition and compute their distinguishing numbers.

	A graph is {\it split} if its vertex set can be partitioned into $A \cup B$ so that the vertices of $A$ induce a complete graph while the vertices of $B$ induce a stable set. The pair $(A,B)$ is referred to as the $KS$-partition of the graph.  Given a split graph $G$ with its $KS$-partition $(A,B)$ and another graph $H$ so that $V(G) \cap V(H) = \emptyset$,  the {\it composition} of $(G,A,B)$ and $H$ is the graph   $(G,A,B) \circ H$ whose vertex set is $V(G) \cup V(H)$ and whose edge set is $E(G) \cup E(H) \cup \{uv | u \in A, v \in V(H) \}$.   A graph that can result from this kind of composition  is {\it decomposable}; otherwise, it is {\it indecomposable}.  Tyshkevich's canonical decomposition theorem \cite{Ty00} states that {\it every} graph can be expressed as  $$(G_r, A_r, B_r) \circ (G_{r-1}, A_{r-1}, B_{r-1}) \circ \cdots \circ (G_1, A_1, B_1) \circ G_0$$
   where  each $G_i$ is an indecomposable graph.  Furthermore, the decomposition is unique up to isomorphism.  
	
	In this paper, we consider a {\it compact} version of Tyshkevich's canonical decomposition where trivial components (i.e., components with only one vertex) are maximally combined to form either a complete graph or a graph of isolated vertices.  We then prove that if the compact canonical decomposition of a graph $G$ is $$(G'_k, A'_k, B'_k) \circ (G'_{k-1}, A'_{k-1}, B'_{k-1}) \circ \cdots \circ (G'_1, A'_1, B'_1) \circ G'_0$$ then a labeling $\phi$ of $G$ is distinguishing if and only if $\phi$ is distinguishing when restricted to $V(G'_i)$ for $i = 0, \hdots, k$.  Consequently, $D(G) = \max \{D(G'_i), i = 0, \hdots, k\}$.  A nice implication of this result is that if $G$ is a threshold graph, then $D(G)$ is equal to size of the largest component in its compact canonical decomposition. 
	
	Next, we apply the above result to unigraphs, utlilizing Tyshkevich's classification of indecomposable unigraphs in \cite{Ty00}, to show that the distinguishing numbers of unigraphs can be computed in linear time.  When $G$ is a unigraph, our algorithm first computes $G$'s compact canonical decomposition,  solves for each component's distinguishing number and then outputs the largest value. 

	We note that the part of our algorithm where it computes the components in $G$'s compact canonical decomposition and identifies their types 
has many similarities to Borri et al.'s linear-time recognition algorithm for unigraphs in \cite{BoCaPe11}.  Their algorithm also goes through the same steps -- albeit for the canonical decomposition and not the compact canonical decomposition -- to determine if the input graph is a unigraph.  Unfortunately, we cannot use their algorithm because our end goals are different\footnote{In their case, the algorithm outputs a red/black edge coloring of the graph when the input is a unigraph and ``failure" otherwise. In our case, the output is the distinguishing number of the input graph.}  so we designed our own algorithm instead.

	Determining a graph parameter's value in a unigraph based on the corresponding values of its indecomposable components is also not unique to our work.  Calamoneri and Petreschi \cite{CaPe11} applied this technique to bound the number of colors needed in an $L(2,1)$-labeling of a unigraph.  They showed that the number of colors they used is at most $1.5 \times OPT$. 
	Nakahata \cite{Na22} also used the technique to show that the clique width of a unigraph is at most $4$.   Perhaps what's surprising is the same idea  did not just provide us with a bound -- it led us to the exact value of the distinguishing number of a unigraph. 
 
         We discuss background information and present preliminary results in section 2.  We prove our first main result -- that a graph's compact canonical decomposition can be used to derive its distinguishing number -- in section 3.   We present  our second main result -- a linear-time algorithm for computing the distinguishing number of a unigraph -- in section 4 and  conclude in section 5.

\section{Background and Preliminary Results}
	
	The {\it degree sequence} $(d_1, d_2, \hdots, d_n)$ of a graph  consists of the degrees of its $n$ vertices ordered from largest to smallest.  A {\it unigraph} $G$ is a graph whose degree sequence is unique up to isomorphism. That is, $G$ and $H$ have the same degree sequences if and only if $G \cong H$.  Unigraphs contain important graph families like threshold graphs \cite{ChHa77, MaPe95}, split matrogenic graphs \cite{HaZv04}, matroidal graphs \cite{Pe77} and matrogenic graphs \cite{FoHa78}.  
Simple examples of unigraphs include complete graphs, $mK_2$ ($m$ copies of $K_2$), $C_5$ and their complements.  We shall encounter unigraphs with more complicated structures in the latter part of the section.  Research on this graph class has been ongoing since the 1970's (e.g. \cite{Li75, Ko76}).  According to Tyskevich \cite{Ty00}, the number of unlabeled unigraphs with $n$ vertices is between $2.3^{n-2}$ and $2.6^n$.

	
	For ease of discussion, we will often abbreviate a degree sequence as $(d_1^{r_1}, \hdots, d_s^{r_2})$ with $d_1  > d_2 > \cdots > d_s$ to indicate that the graph has $r_i$ vertices with degree $d_i$ for $i = 1, \hdots, s$.  Thus, the abbreviated degree sequence of $mK_2$ is $(1^{2m})$.

        A {\it split graph}  is a graph whose vertex set can be partitioned into two sets $A$ and $B$ so that the vertices in $A$ induce a complete graph while the vertices in $B$ induce a stable or independent set.  The pair $(A,B)$ is referred to as a {\it $KS$-partition} of the graph.  Hammer and Simeone \cite{HaSi81} describe the different types of $KS$-partitions that a split graph can have.  
       
\begin{theorem}
(\cite{HaSi81}) For any $KS$-partition $(A,B)$ of a split graph $G$, exactly one of the following holds:
\begin{itemize}
\item $|A| = \omega(G)$ and $|B| = \alpha(G)$
\item $|A| = \omega(G)$ and $|B| = \alpha(G)-1$ (called $K$-max)
\item $|A| = \omega(G) -1$ and $|B| = \alpha(G)$ (called $S$-max)

\end{itemize}
where $\omega(G)$ and $\alpha(G)$ are the clique and independence numbers of $G$ respectively.  Moreover, in a $K$-max partition, there  is a vertex $u \in A$ so that $B \cup \{u\}$ induces a stable set while in an $S$-max partition, there is a vertex $v \in B$ so that $A \cup \{v\}$ induces a complete graph.  \end{theorem}

When a split graph has the first kind of $KS$-partition, it is {\it balanced}; otherwise, it is {\it unbalanced}.  We note below that a balanced split graph has a unique $KS$-partition.   On the other hand, an unbalanced split graph has both a $K$-max and a $S$-max partition because of the existence of the vertices $u$ and $v$. The latter are referred to as {\it swing vertices} of the graph.

\begin{proposition}
\label{uniqueKSprop}
(\cite{ChCoTr16}, Proposition 12) A balanced split graph $G$ has a unique $KS$-partition.  
\end{proposition}

\begin{proposition}
\label{swingvertexprop}
Let $G$ be a split graph with $KS$-partition $(A,B)$.  Then $u \in A$ is a swing vertex if $deg(u) = |A|-1$  and $v \in B$ is a swing vertex if $deg(v) = |A|$. 
\end{proposition}

\begin{proof}
When $u \in A$ has $deg(u) = |A| - 1$, $u$ is not adjacent to any of the vertices in $B$ so $B \cup \{u\}$ is a stable set of $G$. Similarly, when $v \in B$ has $deg(v) = |A|$, $v$ is adjacent  to all the vertices in $A$ so $A \cup \{b\}$ forms a clique in $G$.
\end{proof}

 Throughout our discussion, we use  $(G, A, B)$ to denote a split graph $G$ and its $KS$-partition $(A,B)$.  Its {\it paired degree sequence} is $(d_1,  \hdots, d_p; d_{p+1},  \hdots, d_n)$ where $(d_1,  \hdots, d_p)$ and $(d_{p+1}, \hdots, d_n)$ consists of the degrees of the vertices in $A$ and $B$ respectively, ordered from largest to smallest.  It turns out that we can recognize split graphs from their degree sequences.  

 \begin{theorem}
 \label{splitthm}
 (\cite{HaSi81})  Suppose the degree sequence of graph $G$ is $(d_1, d_2, \hdots, d_n)$ with $n \ge 2$ and $h = \max\{i: d_i \ge i-1\}$.  Then $G$ is a split graph if and only if 
  $$ \sum_{i=1}^h d_i = h(h-1) + \sum_{i = h+1}^n d_i.$$  If the equality above holds, then $h = \omega(G) $.  Furthermore, a $KS$-partition of the graph  has the paired degree sequence  $(d_1, d_2, \hdots, d_h; d_{h+1}, \hdots, d_n)$.  
 \end{theorem}       
   
   In the appendix, we present the algorithm DetermineSplit based on Theorem \ref{splitthm}.  It takes $G$'s degree sequence as input and outputs the paired degree sequence  $(d_1, d_2, \hdots, d_h; d_{h+1}, \hdots, d_n)$  if and only if $G$ is a split graph.  
   
   \begin{corollary}
   \label{splitalgcor}
   Given the degree sequence of a graph with $n$ vertices, DetermineSplit determines if the graph is split in $O(n)$ time.
   \end{corollary}
   
\subsection{Automorphisms and Distinguishing Numbers}
   
Recall that an {\it automorphism} of $G$ is a bijection $\pi:V(G) \rightarrow V(G)$ that preserves the adjacencies of $G$. 
We say that $\pi$ is {\it trivial} if $\pi$ is the identity map; otherwise, $\pi$ is {\it non-trivial}.  The automorphism $\pi$ {\it fixes} a subset $V' \subseteq V(G)$ if and only if for every $v \in V'$, $\pi(v) \in V'$.  When $G$ is a split graph with $KS$-partition $(A,B)$,  every automorphism of $(G,A,B)$ must fix $A$ and $B$ too.  Let  $Aut(G)$ and $Aut(G,A,B)$ contain the automorphisms of $G$ and $(G,A,B)$ respectively. 
    It is obvious that $Aut(G,A,B) \subseteq Aut(G)$.  In the next proposition, we describe a situation when they are equal to each other.
   
   \begin{proposition}
   \label{equalautomorph}
   When $G$ is a balanced split graph with $KS$-partition $(A,B)$, $Aut(G) = Aut(G,A,B)$.
   \end{proposition}
   
    \begin{proof}
    By Propositions \ref{uniqueKSprop},  $(A,B)$ is the unique $KS$-partition of $G$.  Thus, every automorphism of $G$ must fix $A$ and $B$.  That is, $Aut(G) \subseteq Aut(G, A, B)$.  Since it is also the case that $Aut(G,A,B) \subseteq Aut(G)$, the proposition follows.    \end{proof}

Let $\phi: V(G) \rightarrow \{1, 2, \hdots, c\}$ be a $c$-labeling of $G$ that uses $c$ colors.  We say that $\phi$ is a {\it distinguishing labeling} if,  for any non-trivial automorphism $\pi$ of $G$,  there is some vertex $v$ so that $\phi(v) \neq \phi(\pi(v))$.  That is, the only automorphism of the labeled graph $(G, \phi)$ is the identity map.  The {\it distinguishing number of $G$}, $D(G)$, is the smallest $c$ such that $G$ has a distinguishing $c$-labeling.   For example, $D(K_n) = n$ because no two vertices of a complete graph can be assigned the same color in a distinguishing labeling.  On the other hand,  $D(P_n) = 2$ for $n \ge 2$ because it suffices to assign distinct colors to the endpoints of the path to distinguish its vertices.  When $G$ is a split graph with $(A,B)$ as a $KS$-partition, we extend the notion of  distinguishing labelings and distinguishing number of $(G,A,B)$ with respect to the automorphisms of $(G,A,B)$.   As an example, consider the graph $G = K_3 \cup K_1$.  If $A$ consists of the three vertices of $K_3$ and $B$ the vertex of $K_1$, then $D(G,A,B) = 3$, which is also $D(G)$.   But if $A$ consists of two vertices of $K_3$ and $B$ consists of the remaining vertex of $K_3$ and the vertex of $K_1$, then $D(G,A,B)  = 2$ because we only have to assign distinct colors to the two vertices in $A$.

\begin{proposition}
\label{samedistprop}
 Let $G$ and $H$ be two graphs on the same set of vertices; i.e., $V(G) = V(H)$. If $Aut(G) = Aut(H)$ then $D(G) = D(H)$. 
\end{proposition}

\begin{proof}
Let $\phi$ be a distinguishing labeling of $G$.  Suppose $\phi$ is not distinguishing for $H$.  Then there is some non-trivial automorphism $\pi$ of $H$ such that $\phi(v) = \phi(\pi(v))$ for each $v \in V(H)$.  But $\pi \in Aut(G)$ too so $\phi$ cannot be a distinguishing labeling of $G$, a contradiction.  The same argument shows that a distinguishing labeling for $H$ is also one for $G$. It follows that $D(G) = D(H)$. 
\end{proof}


Combining Proposition \ref{equalautomorph} with the argument in the proof of Proposition \ref{samedistprop}, the following is also true. 

\begin{proposition}
\label{samedistprop2}
 When $G$ is a balanced split graph with $KS$-partition $(A,B)$, $D(G) = D((G,A,B))$. 
\end{proposition}

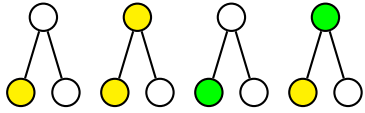
\begin{figure}

 \begin{center}
\begin{tikzpicture}

\begin{scope}[scale=.5,auto=center,every node/.style={circle,thick,draw,minimum size=0.05cm}] 
    \node (a1) at (0, 0) {};
    \node[fill = yellow] (a2) at (-0.6,-2) {};
    \node   (a3) at (0.6,-2) {}; 
     \node[fill = yellow] (b1) at (2.5, 0) {};
    \node[fill = yellow] (b2) at (1.9,-2) {};
    \node (b3) at (3.1, -2) {};
     \node (c1) at (5, 0) {};
    \node[fill=green] (c2) at (4.4,-2) {};
    \node (c3) at (5.6,-2) {}; 
     \node[fill=green] (d1) at (7.5, 0) {};
    \node[fill=yellow] (d2) at (6.9,-2) {};
    \node (d3) at (8.1, -2) {};
\end{scope}
\path[-, draw, thick] (a1) -- (a2);
\path[-, draw, thick] (a1) -- (a3);
\path[-, draw, thick] (b1) -- (b2);
\path[-, draw, thick] (b1) -- (b3);
\path[-, draw, thick] (c1) -- (c2);
\path[-, draw, thick] (c1) -- (c3);
\path[-, draw, thick] (d1) -- (d2);
\path[-, draw, thick] (d1) -- (d3);
\end{tikzpicture}
\end{center}
\caption{Four inequivalent distinguishing labelings of $K_{1,2}$ using three colors.}
\label{4distlabelfig}
\end{figure}

When $G = mH$ (i.e., $G$ is made up of $m$ copies of graph $H$), it is useful to know $D(H)$ but it is not enough to determine $D(G)$.   Instead, we have to consider the number of inequivalent distinguishing labelings of $H$ -- a concept used in  \cite{ArDe04} and \cite{Ch06}  independently to compute the distinguishing numbers of trees and forests.    Let $\phi$ and $\phi'$ be two distinguishing labelings of $G$.  We say that  $\phi$ and $\phi'$ are {\it equivalent} if $(G, \phi) \cong (G, \phi')$; otherwise, they are {\it inequivalent}.  Figure \ref{4distlabelfig} shows  four inequivalent distinguishing labelings of the star $K_{1, 2}$.  Given $c$ colors, define $D(G, c)$ as the number of inequivalent distinguishing $c$-labelings of $G$.

\begin{proposition}
Let $G = mH$.  Then $D(G) = \min\{c: D(H,c) \ge m\}$.   
\end{proposition}


Later,  we will determine the distinguishing numbers of unigraphs.  Two graphs play an important role: $mK_2$ and $S(p,q)$.  We show how to compute their distinguishing numbers below.
\medskip


\noindent {\it Example 1:}  Consider $mK_2$.  Every distinguishing labeling of $K_2$ must assign distinct colors  to its two vertices.  Furthermore, two distinguishing labelings of $K_2$ are inequivalent if they assign different pairs of colors to the vertices of $K_2$.  
Thus, $D(K_2, c) = \binom{c}{2}$ and $D(mK_2) = \min\{c: \binom{c}{2} \ge m\}$. 
\medskip

\noindent {\it Example 2:}  Let $S(p,q)$ denote the graph that is made up of $q \ge 2$ copies of the star $K_{1, p}$ whose centers are pairwise adjacent to each other.  See Figure \ref{TyshFig1} for an example. 
 Let $K^*_{1,p}$ be $K_{1,p}$ whose center is designated as a root so that every automorphism of $K^*_{1,p}$ maps the center to itself.\footnote{This extra step of defining $K^*_{1,p}$ is unnecessary when $p \ge 2$ because $K^*_{1,p} = K_{1,p}$.  We only defined it to deal with the case when $p = 1$ so that we can designate one vertex as the ``center" of the star and the other vertex as a leaf of the star.  Otherwise, if we just considered $K_{1,1}$, both of its vertices can play the role of the center of the star.} 
 By definition, $qK^*_{1,p}$ is a subgraph of $S(p,q)$, and it is straightforward to verify that  $Aut(S(p,q)) = Aut(q K^*_{1,p})$.     By Proposition \ref{samedistprop},  $D(S(p,q)) = D(q K^*_{1,p})$.  To distinguish $K^*_{1,p}$, it is enough that distinct colors be assigned to the leaves of $K^*_{1,p}$. Furthermore, two distinguishing labelings $\phi$ and $\phi'$ of $K^*_{1,p}$  are inequivalent if and only if $\phi$ and $\phi'$ assign different colors to the center of $K^*_{1,p}$ or they assign distinct sets of colors to the leaves of $K^*_{1,p}$.  Thus, $D(K^*_{1,p}, c) = c \binom{c}{p}$ so $D(S(p,q)) =  \min\{c: c \binom{c}{p} \ge q \}.$


\subsection{Tyshkevich's canonical decomposition and its compact version}


   Given a split graph $(G,A,B)$ and another graph $H$ so that $V(G) \cap V(H) = \emptyset$,  Tyshkevich \cite{Ty00} defined the {\it composition} of $(G,A,B)$ and $H$ as the graph  $(G,A,B) \circ H$ whose vertex set is $V(G) \cup V(H)$ and whose edge set is $E(G) \cup E(H) \cup \{uv | u \in A, v \in V(H) \}$.   If, additionally, $H$ has a $KS$-partition $(C,D)$, then $(G, A, B) \circ (H, C, D) = (F, A\cup C, B \cup D)$ where $F = (G, A,B) \circ H$. The composition can be applied multiple times to produce a new graph; e.g., 
   $(G_r, A_r, B_r) \circ (G_{r-1}, A_{r-1}, B_{r-1}) \circ \cdots \circ (G_1, A_1, B_1) \circ G_0.$
 \noindent Based on the definition, our inclination is to apply the composition operation from right to left.  But in fact, the composition operation is associative; it can be applied in any order.

   Let us refer to a graph as {\it decomposable} if  there is some split graph $(G,A,B)$ and another graph $H$,  both with at least one vertex each,  so that the graph is equal to $(G,A,B) \circ H$.  Otherwise, the graph is {\it indecomposable}.

 A vertex is {\it isolated} in $G$ if it is not adjacent to any vertices; it is {\it dominant} if it is adjacent to all the vertices except itself.  When $G$ is made up of just one vertex $v$, then we simply refer to it as $\{v\}$.  The following results are well known; we prove part c for completeness.

   \begin{proposition}
   \label{decomposable}  Let $G$ be a graph with two or more vertices. 
\begin{enumerate}
   
\item[a.]  If $G$ has an isolated vertex $u$ then $G = (\{u\}, \emptyset, \{u\}) \circ (G - u)$. 
  
\item[b.]  If $G$ has a dominant vertex $v$  then $G = (\{v\}, \{v\}, \emptyset) \circ (G-v)$. 
   
 \item[c.]  If $G$ is an indecomposable split graph then $G$ is balanced.
 
  \end{enumerate}
    
 
   \end{proposition}
   
  \begin{proof}
  Parts a and b are obvious so let us prove part c.  Suppose $G$ is an indecomposable split graph but $G$ is unbalanced.  Then it has a $KS$-partition $(A,B)$ that is $K$-max.  Moreover, some $w \in A$ is a swing vertex so $w$ is adjacent to all the vertices in $A$ but not to any vertices in $B$.  Thus, $G = (G-w, A-\{w\}, B) \circ \{w\}$, contradicting the fact that $G$ is indecomposable.  Hence, $G$ has to be balanced. 
   \end{proof} 
   
   
 
 \medskip
   
   Tyshkevich's result states why indecomposable graphs are important.

   \begin{theorem}
   \label{canon}
   (The Canonical Decomposition Theorem of Tyshkevich \cite{Ty00})  Every graph $G$ can be expressed as 
   $$G = (G_r, A_r, B_r) \circ (G_{r-1}, A_{r-1}, B_{r-1}) \circ \cdots \circ (G_1, A_1, B_1) \circ G_0$$
   where  each $G_i$ is an indecomposable graph.  Furthermore, the decomposition is unique up to isomorphism.  That is, let $$G' = (G'_\ell, A'_\ell, B'_\ell) \circ (G'_{\ell-1}, A'_{\ell-1}, B'_{\ell-1}) \circ \cdots \circ (G'_1, A'_1, B'_1) \circ G'_0.$$ Then $G \cong G'$ if and only if $r = \ell$, $G_0 \cong G'_0$ and $(G_i, A_i, B_i) \cong (G'_i, A'_i, B'_i)$ for $i = 1, \hdots, r$. 
   
   \end{theorem}

Every component $G_i$, $i > 0$, in Tyshkevich's decomposition is an indecomposable split graph and has a unique $KS$-partition.  Thus, it is common to just state the decomposition of $G$ as $G_r \circ G_{r-1} \circ \cdots \circ G_1 \circ G_0$. 
In her paper (\cite{Ty00}, paragraph 1),  Tyshkevich hinted that the canonical decomposition of a graph can be computed in linear time.  We formally prove this fact in the appendix.  It is also implied in the work of Borri et al. \cite{BoCaPe11}. 

\begin{theorem}
\label{decomposealgthm}
Let $G$ be a graph with $n$ vertices and $m$ edges.  Suppose the canonical decomposition of $G$ is  $G = (G_r, A_r, B_r) \circ (G_{r-1}, A_{r-1}, B_{r-1}) \circ \cdots \circ (G_1, A_1, B_1) \circ G_0.$. The algorithm DECOMPOSE($G$) returns a stack $S$ that contains the (paired) degree sequences of the $G_i$'s in order in $O(n+m)$ time.

\end{theorem}

Call a component $G_i$, $i \ge 0$,  {\it trivial}  if it has only one vertex; otherwise, it is {\it non-trivial}.  There are two types of trivial split graph components depending on whether the single vertex is in the $K$-part or $S$-part of the $KS$ partition. Let $K_1$ and $S_1$ refer to the two types respectively whose paired degree sequences are $(0; \emptyset)$ and $(\emptyset; 0)$. 
 When $G_0$ is also a trivial component (i.e., its degree sequence is $(0)$), we say its type is $K_1$ if $G_1$ has type $K_1$ and $S_1$ if $G_1$ has type $S_1$.   That is, $G_0$ follows $G_1$'s type if both of them are trivial components.  Otherwise, when $G_1$ is  a non-trivial component or does not exist, we leave $G_0$ alone. It is possible that consecutive components of $G$ are trivial and of the same type.  We want to combine them {\it maximally}.

\begin{definition}
\label{compactdef}
The {\it compact canonical decomposition of $G$} is formed by maximally combining consecutive trivial components of the same type in the canonical decomposition of $G$.  That is, if there are $m$ consecutive trivial components of type $K_1$, replace them with the complete graph on $m$ vertices with $K$-max as its $KS$-partition (i.e., the degree sequence is $((m-1)^m; \emptyset)$).  On the other hand, if their type is $S_1$, replace them with the graph with $m$ isolated vertices whose $KS$-partition is $S$-max (i.e., the degree sequence is $(\emptyset; 0^m)$).
\end{definition}

\medskip

\noindent {\it Example 3:}  Consider $G$ whose canonical decomposition is made up entirely of trivial components:
$$G = (0; \emptyset) \circ (0; \emptyset) \circ (\emptyset; 0) \circ (\emptyset; 0) \circ (\emptyset; 0) \circ (0).$$
By combining the trivial components of the same type (and where $G_0$ adopts the type of $G_1$), we obtain the compact canonical decomposition of $G$:
$$G = (1^2; \emptyset) \circ (\emptyset; 0^4).$$

\begin{corollary}
 \label{canon2}
(The Compact Canonical Decomposition Theorem)    The compact canonical decomposition of a graph $G$ is unique up to isomorphism where each component is an indecomposable graph with at least two vertices, a complete graph or a graph of isolated vertices. 
\end{corollary}

\begin{proof} 
The compact canonical decomposition of $G$ is derived from Tyshkevich's canonical decomposition by combining consecutive trivial components of the same type maximally. Thus, a component of Tyshkevich's decomposition is either untouched or combined with other trivial components.  In the former case, the component is an indecomposable graph with at least two vertices or a single node.  In the latter case, a complete graph or a graph of isolated vertices is formed.   Finally, since there is only one way to do the maximal combination step and Tyshkevich's decomposition is unique up to isomorphism, it follows that the compact canonical decomposition of $G$ is also unique up to isomorphism.  \end{proof}



In the appendix, we present DECOMPOSE-C($G$) which takes the output of DECOMPOSE($G$) and maximally combines the trivial components of the same type to produce the compact canonical decomposition of $G$.  Like DECOMPOSE($G$), the algorithm runs in linear time. 

\begin{theorem}
\label{decompose-c}
Let $G$ be a graph with $n$ vertices and $m$ edges.  Suppose the compact canonical decomposition of $G$ is  $G = (G_k, A_k, B_k) \circ (G_{k-1}, A_{k-1}, B_{k-1}) \circ \cdots \circ (G_1, A_1, B_1) \circ G_0.$
The algorithm DECOMPOSE-C($G$) returns a stack $T$ that contains the abbreviated (paired) degree sequences of $G$ in order in $O(n+m)$ time.


\end{theorem}

\subsection{Indecomposable Unigraphs and their Distinguishing Numbers}

 The following corollary of Tyshkevich's canonical decomposition theorem is perhaps not surprising.
      
 \begin{corollary} (\cite{Ty00}) $G$ is a unigraph if and only if each indecomposable component in the canonical decomposition of $G$ is also a unigraph.
 
 \end{corollary}      
 

A complete graph and a graph of isolated vertices are unigraphs so the following is true as well.
 
  \begin{corollary}
  $G$ is a unigraph if and only if each component in the compact canonical decomposition of $G$ is also a unigraph.
 
  \end{corollary}      

 More importantly for us, Tyshkevich identified all the indecomposable unigraphs. But before we present them, let us consider two graph operations.  Recall that $\overline{G}$, the {\it complement} of graph $G$, has the same vertex as $G$ and two vertices are adjacent if and only if they are not adjacent in $G$.  When $G$ is a split graph with $KS$-partition $(A,B)$,  the complement of $(G,A,B)$ is $\overline{(G,A,B)} = (\overline{G}, B, A)$. Its {\it inverse} is $(G,A,B)^I = (G^I, B, A)$ where $G^I$ is obtained from $G$ by deleting the edges $\{aa': a, a' \in A\}$ and adding the edges $\{bb': b, b' \in B\}$.

 For the rest of this paper, we will primarily be dealing with indecomposable graphs.   When such a graph $G$ is split, it is balanced and therefore has only one $KS$-partition.  To simplify the discussion, we will often refer its complement and inverse as $\overline{G}$ and $G^I$ respectively.  It is easy to verify that $(G^I)^I = G$ and $\overline{G^I} = \overline{G}^I$. 
 

 \begin{proposition}
  If graph $G$ is indecomposable, then so is its complement $\overline{G}$.  If $G$ is an indecomposable split graph then so is $G^I$. 
   \end{proposition}
   

\begin{proposition}
For any graph $G$, $D(G) = D(\overline{G})$.  Additionally, when $G$ is an indecomposable split graph, $D(G) = D(G^I)$.
\label{equaldistnums}
\end{proposition}
 
 \begin{proof}
 It is obvious that $Aut(\overline{G}) = Aut(G)$.  By Proposition \ref{samedistprop}, $D(\overline{G}) = D(G)$.  To prove the second part, we will show that $Aut(G^I) = Aut(G)$.   Let  $\pi \in Aut(G)$ and $(A,B)$ be the unique $KS$-partition of $G$.
 
First, assume $u \in A$ and $v \in B$.  We know that $uv$ is an edge of $G$ if and only if $\pi(u) \pi(v)$ is an edge of $G$.  Since $G$ and $G^I$ have exactly the same set of edges between $A$ and $B$, it follows that $uv$ is an edge of $G^I$ if and only if $\pi(u)\pi(v)$ is an edge of $G^I$.  

Next, let $u, v \in A$ or $u,v \in B$.  By Proposition \ref{equalautomorph}, $\pi$ fixes $A$ and $B$ so $\pi(u), \pi(v) \in A$ in the former case while $\pi(u), \pi(v) \in B$ in the latter case.  In $G^I$,  no pair of vertices in $A$ are adjacent while every pair of vertices in $B$ are adjacent.  Thus, for the two cases, $\pi$ preserves the adjacencies of $G^I$.  We have shown that $\pi \in Aut(G^I)$ so $Aut(G) \subseteq Aut(G^I)$. But $(G^I)^I = G$. The argument above also implies that $Aut(G^I) \subseteq Aut(G)$ so $Aut(G) = Aut(G^I)$.  By Proposition \ref{samedistprop}, $D(G^I) = D(G)$. \end{proof}

\begin{figure}

 \begin{center}
\begin{tikzpicture}

\begin{scope}[scale=.5,auto=center,every node/.style={circle,thick,draw,minimum size=0.05cm}] 
    \node (a1) at (0, 0) {};
    \node (a2) at (0,-2) {};
     \node (b1) at (1, 0) {};
    \node (b2) at (1,-2) {};
    \node(c1) at (3.5,0) {};
    \node(c2) at (2.5,-2) {};
    \node(c3) at (3.5,-2) {};
    \node(c4) at (4.5,-2) {};
\end{scope}
\node(l1) at (1, -2) {$U_2(2,3)$};
\path[-, draw, thick] (a1) -- (a2);
\path[-, draw, thick] (b1) -- (b2);
\path[-, draw, thick] (c1) -- (c2);
\path[-, draw, thick] (c1) -- (c3);
\path[-, draw, thick] (c1) -- (c4);
\end{tikzpicture}
\hspace*{0.3in}
\begin{tikzpicture}
\begin{scope}[scale=.5,auto=center,every node/.style={circle,thick,draw,minimum size=0.05cm}] 
    \node (a1) at (0, 1) {};
    \node (a2) at (1.25,2) {};
     \node (a3) at (1.25,0) {};
      \node (a4) at (2.5,1) {};
     \node(b1) at (3.3, 2.75) {}; 
     \node(b2) at (4.2, 1.5){};
     \node(c1) at (4.2, 0.5){};
     \node(c2) at (3.3, -0.5){};
\end{scope}
\node(l1) at (1, -1) {$U_3(2)$};
\path[-, draw, thick] (a1) -- (a2);
\path[-, draw, thick] (a1) -- (a3);
\path[-, draw, thick] (a4) -- (a2);
\path[-, draw, thick] (a4) -- (a3);
\path[-, draw, thick] (a4) -- (b1);
\path[-, draw, thick] (a4) -- (b2);
\path[-, draw, thick] (a4) -- (c1);
\path[-, draw, thick] (a4) -- (c2);
\path[-, draw, thick] (b1) -- (b2);
\path[-, draw, thick] (c1) -- (c2);
\end{tikzpicture}
\hspace*{0.4in}
\begin{tikzpicture}
\begin{scope}[scale=.5,auto=center,every node/.style={circle,thick,fill=gray,draw,minimum size=0.05cm}] 
  \node (a1) at (0, 0) {};
   \node (b1) at (1.5, 0) {};
   \node(c1) at (3,0) {};
 \end{scope}
\begin{scope}[scale=.5,auto=center,every node/.style={circle,thick,draw,minimum size=0.05cm}] 
    \node (a2) at (0,-2) {};
    \node (b2) at (1.5,-2) {};
    \node(c2) at (3,-2) {};    
\end{scope}
\node(l1) at (0.8, -2) {$S(1, 3)$};
\path[-, draw, thick] (a1) -- (a2);
\path[-, draw, thick] (b1) -- (b2);
\path[-, draw, thick] (c1) -- (c2);
\path[-, draw, thick] (a1) -- (b1);
\path[-, draw, thick] (b1) -- (c1);
\draw[-, draw, thick]  (a1) to[out=45, in=135] (c1);
\end{tikzpicture}
\end{center}
\vspace*{-1em}
\caption{We show examples of the indecomposable unigraphs $U_2(m, \ell)$, $U_3(m)$ and $S(p,q)$. The gray nodes of $S(1, 3)$ are the ones that are part of the clique in the $KS$-partition of the graph.}
\label{TyshFig1}
\end{figure}
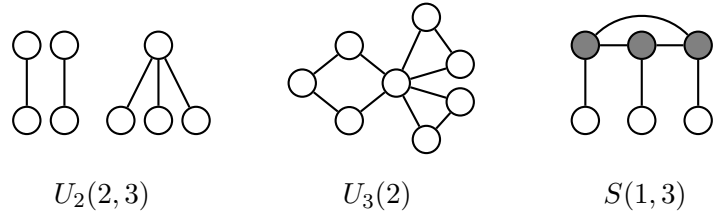

   We start with the indecomposable unigraphs that are not split graphs.  Table \ref{UnsplittableTable} lists the graphs mentioned by Tyshkevich (Theorem \ref{Tyshkevich-unsplittable}), their degree sequences and distinguishing numbers.  Aside from the fact that $mK_2, m \ge 2$, is itself an indecomposable non-split unigraph, it also serves as a building block for the other graphs.  The graph $U_2(m, \ell)$ is the disjoint union of $mK_2$ and the star $K_{1, \ell}$ while in $U_3(m)$, there is a center vertex that is adjacent to the endpoints of all the vertices of $mK_2$  to form $m$ triangles and to the endpoints of a $3$-vertex path to form a $4$-cycle. (See  Figure \ref{TyshFig1} for some examples.) Thus, we express $D(U_2(m, \ell)$ and $D(U_3(m))$ in terms of $D(mK_2)$.   
 \begin{theorem}
 \label{Tyshkevich-unsplittable}
 (Tyshkevich \cite{Ty00}) Let $G$ be an indecomposable non-split unigraph with two or more vertices.  Then $G$ or $\overline{G}$ is one of the following graphs: 
 $C_5, m K_2$ with $m \ge 2$, $ U_2(m,\ell)$ with $m \ge 1, \ell \ge 2$ or $U_3(m)$ with $m \ge 1$. 
 \end{theorem}
 
 \begin{lemma}
 The degree sequences and distinguishing numbers of the graphs in Table \ref{UnsplittableTable} are correct. 
  \end{lemma}
 
 \begin{proof}  It is straightforward to check that the degree sequences listed in the table are correct. We will just verify the distinguishing numbers of the graphs.   Albertson and Collins  \cite{AlCo96} noted that $D(C_5) = 3$ in their introductory paper on distinguishing labelings.  In Section 2.1, we showed that $D(m K_2) = \min \{ c: \binom{c}{2} \ge m \}$. For $U_2(m, \ell)$, since $\ell \ge 2$,  no automorphism of the graph will map the vertices of $mK_2$ to the vertices of $K_{1,\ell}$.  Thus, $D(U_2(m,\ell)) = \max \{D(mK_2), D(K_{1,\ell})\}$. We already described how to compute $D(mK_2)$.  Since the leaves of $K_{1,\ell}$ have to be assigned distinct colors, $D(K_{1,\ell}) = \ell$.  
    
    Finally, in $U_3(m)$, let $v$ be the center vertex. Every automorphism of $U_3(m)$ fixes $v$ since it has the largest degree in the graph.  Thus,  every distinguishing labeling of $U_3(m)$ has to  distinguish $mK_2$ and the $3$-vertex path connected connected to $v$; i.e., $D(U_3(m)) = \max \{D(mK_2), D(P_3) \}$.   But $D(P_3) = 2$ and $D(m K_2) \ge 2$ so $D(U_3(m)) = D(m K_2)$.  \end{proof}

 \begin{table}
 \begin{center}
\small
 \begin{tabular}{|m{2.6in}|m{1.5 in}|m{1.2in}|}
 \hline
{Indecomposable non-split graph $G$} & Degree sequence of $G$ & ${D(G)}$  \\
\hline
$C_5$ & $(2^5)$  & 3 \\
\hline
$mK_2$ with $m \ge 2$ & $(1^{2m})$ &  $\min \{c: \binom{c}{2} \ge m \}$ \\
\hline
$U_2(m,\ell) = mK_2 \cup K_{1,\ell}$ with $m \ge 1, \ell \ge 2$
& $(\ell, 1^{2m + \ell})$  & 
$\max\{D(mK_2), \ell \}$ \\
\hline
$U_3(m)$, $m \ge 1$ &  $(2m+2, 2^{2m+3})$
& $D(mK_2)$  \\
\hline
\end{tabular} 
 \caption{The four types of indecomposable non-split unigraphs, their degree sequences and distinguishing numbers.}
 \label{UnsplittableTable}
 \end{center}
 \end{table}

\begin{figure}[t]

 \begin{center}
\begin{tikzpicture}
\begin{scope}[scale=.5,auto=center,every node/.style={circle,thick,fill=gray,draw,minimum size=0.05cm}] 
  \node (a1) at (0, 0) {};
   \node (b1) at (1.5, 0) {};
   \node(c1) at (3,0) {};
   \node(e1) at (6,0) {};
   \node(f1) at (8,0){};
 \end{scope}
 \begin{scope}[scale=.5,auto=center,every node/.style={circle,thick,fill=yellow,draw,minimum size=0.05cm}] 
  \node(d2) at (4.5, -2){};
 \end{scope}
\begin{scope}[scale=.5,auto=center,every node/.style={circle,thick,draw,minimum size=0.05cm}] 
    \node (a2) at (0,-2) {};
    \node (b2) at (1.5,-2) {};
    \node(c2) at (3,-2) {};    
    \node(e2) at (5.5, -2){};
    \node(e3) at (6.5, -2){}; 
    \node(f2) at (7.5, -2){};
    \node(f3) at (8.5, -2){};
\end{scope}
  \node(l1) at (2, -2) {$S_3(1, 3, 2)$};
\path[-, draw, thick] (a1) -- (a2);
\path[-, draw, thick] (b1) -- (b2);
\path[-, draw, thick] (c1) -- (c2);
\path[-, draw, thick] (a1) -- (b1);
\path[-, draw, thick] (b1) -- (c1);
\draw[-, draw, thick]  (a1) to[out=30, in=150] (c1);
\path[-, draw, thick] (a1) -- (d2);
\path[-, draw, thick] (b1) -- (d2);
\path[-, draw, thick] (c1) -- (d2);
\path[-, draw, thick] (e1) -- (e2);
\path[-, draw, thick] (e1) -- (e3);
\path[-, draw, thick] (f1) -- (f2);
\path[-, draw, thick] (f1) -- (f3);
\path[-, draw, thick] (c1) -- (e1);
\path[-, draw, thick] (e1) -- (f1);
\draw[-, draw, thick]  (b1) to[out=30, in=150] (e1);
\draw[-, draw, thick]  (b1) to[out=30, in=150] (f1);
\draw[-, draw, thick]  (c1) to[out=30, in=150] (f1);
\draw[-, draw, thick]  (a1) to[out=30, in=150] (e1);
\draw[-, draw, thick]  (a1) to[out=30, in=150] (f1);
\end{tikzpicture}
\hspace*{0.3in}
\begin{tikzpicture}
\begin{scope}[scale=.5,auto=center,every node/.style={circle,thick,fill=gray,draw,minimum size=0.05cm}] 
  \node (a1) at (0, 0) {};
   \node (b1) at (2, 0) {};
   \node(e1) at (6,0) {};
   \node(f1) at (8,0){};
 \end{scope}
 \begin{scope}[scale=.5,auto=center,every node/.style={circle,thick,fill=yellow,draw,minimum size=0.05cm}] 
  \node(d2) at (1, -2){};
 \end{scope}
  \begin{scope}[scale=.5,auto=center,every node/.style={circle,thick,fill=green,draw,minimum size=0.05cm}] 
  \node(c1) at (4, 0){};
 \end{scope}
\begin{scope}[scale=.5,auto=center,every node/.style={circle,thick,draw,minimum size=0.05cm}] 
    \node (a2) at (0,-2) {};
    \node (b2) at (2,-2) {};
    \node(e2) at (5.5, -2){};
    \node(e3) at (6.5, -2){}; 
    \node(f2) at (7.5, -2){};
    \node(f3) at (8.5, -2){};
\end{scope}
  \node(l1) at (2, -2) {$S_4(1,  2)$};
\path[-, draw, thick] (a1) -- (a2);
\path[-, draw, thick] (b1) -- (b2);
\path[-, draw, thick] (c1) -- (a2);
\path[-, draw, thick] (c1) -- (b2);
\path[-, draw, thick] (a1) -- (b1);
\path[-, draw, thick] (b1) -- (c1);
\draw[-, draw, thick]  (a1) to[out=30, in=150] (c1);
\path[-, draw, thick] (a1) -- (d2);
\path[-, draw, thick] (b1) -- (d2);
\path[-, draw, thick] (c1) -- (e2);
\path[-, draw, thick] (c1) -- (f2);
\path[-, draw, thick] (c1) -- (e3);
\path[-, draw, thick] (c1) -- (f3);
\path[-, draw, thick] (e1) -- (e2);
\path[-, draw, thick] (e1) -- (e3);
\path[-, draw, thick] (f1) -- (f2);
\path[-, draw, thick] (f1) -- (f3);
\path[-, draw, thick] (c1) -- (e1);
\path[-, draw, thick] (e1) -- (f1);
\draw[-, draw, thick]  (b1) to[out=30, in=150] (e1);
\draw[-, draw, thick]  (b1) to[out=30, in=150] (f1);
\draw[-, draw, thick]  (c1) to[out=30, in=150] (f1);
\draw[-, draw, thick]  (a1) to[out=30, in=150] (e1);
\draw[-, draw, thick]  (a1) to[out=30, in=150] (f1);
\end{tikzpicture}

\end{center}
\vspace*{-1em}
\caption{We show examples of $S_3(p, q_1, q_2)$ and $S_4(p, q)$.  The vertices $e$ and $f$ are colored yellow and green respectively. }
\label{TyshFig2}
\end{figure}
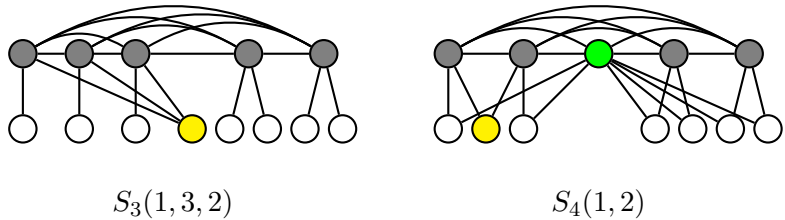

 Next, we consider the indecomposable unigraphs that are also split graphs.  Table \ref{SplittableTable} lists the five graphs mentioned by Tyshkevich (Theorem \ref{Tyshkevich-splittable})  together with their paired degree sequences and distinguishing numbers.  The graph that plays an important role is $S(p,q)$.  As we already noted,  it is made up of $q$ copies of the  star $K_{1,p}$  and whose centers  induce a complete graph.   The graph  $S_2(p_1, q_1, p_2,  q_2, \cdots, p_m, q_m)$ is made up of $S(p_1, q_1)$, $S(p_2, q_2), \hdots, S(p_m, q_m)$ and the centers of the stars from different $S(p_i, q_i)$'s are also pairwise adjacent to each other.   The graph  $S_3(p, q_1, q_2)$, is made up of $S(p, q_1)$ and $S(p+1, q_2)$ and all their stars' centers are pairwise adjacent to each other.  Additionally,  there is a vertex $e$ that is adjacent to the centers of the stars in $S(p, q_1)$ only.  Finally, $S_4(p, q)$  is  made up of $S_3(p, 2, q)$ and a  vertex $f$ that is adjacent to all the vertices of $S_3(p,2, q)$ except for $e$. Figure \ref{TyshFig1} shows $S(1,3)$ while Figure \ref{TyshFig2} shows $S_3(1,3,2)$ and $S_4(1,2)$.

 \begin{theorem}
 \label{Tyshkevich-splittable}
 (Tyshkevich \cite{Ty00}) Let $G$ be an indecomposable split unigraph.  Then $G$,  $\overline{G}$, $G^I$ or $\overline{G^I}$ is one of the following graphs: 
a single-vertex graph,  $S(p,q)$ with $p \ge 1, q \ge 2, $  $ S_2(p_1, q_1, p_2, q_2, \cdots, p_m, q_m)$ with $m \ge 2$ and $p_1 > p_2 > \hdots > p_m,$   $ S_3(p, q_1, q_2)$ with $p \ge 1, q_1 \ge 2, q_2 \ge 1, \;$  and $\;S_4(p,q)$ with $p \ge 1, q \ge 1$.
 
 \end{theorem}
 
 \begin{table}
 \small
\label{Splittable}
 \begin{tabular}{|m{1.7in}|m{2.2in}|m{2in}|}
 \hline
Indecomposable split graph $G$ & Paired degree sequence of $G$ & ${D(G)}$ \\
\hline
Single-vertex graph & $(0; \emptyset )$ or $(\emptyset;0)$ & $1$\\
\hline
$S(p,q)$, $p \ge 1$, $q \ge 2$ & $((p+q-1)^q; 1^{pq})$  & $\min \{c: c\binom{c}{p} \ge q \}$  \\
\hline
$S_2(p_1, q_1, \hdots, p_m, q_m)$,  \hspace*{0.03in} $m \ge 2$ and $p_1 > p_2 > \hdots > p_m$ & $((p_1+ N -1)^{q_1}, \cdots, $ $(p_m + N-1)^{q_m};$ $1^{p_1q_1 + \cdots + p_m q_m})$ where $N = \sum_{i=1}^m q_i$   & $\max\{D(S(p_i, q_i)), i = 1,\hdots,m\}$ \\
\hline
$S_3(p, q_1, q_2)$ with $p \ge 1, q_1 \ge 2, q_2 \ge 1$
& $((p + q_1 + q_2)^{q_1 + q_2}; q_1,$ $1^{pq_1 + (p+1)q_2})$ & $
\max\{D(S(p, q_1)), D(S(p+1, q_2)) \}$ \\
\hline
$S_4(p,q)$ with $p \ge 1$, $q \ge 1$ & $(2(p+q+1)+qp, (p+q+3)^{q+2}; 2^{qp + 2p + q + 1})$ & $\max \{D(S(p,2)), D(S(p+1, q)) \}$  \\
\hline
\end{tabular} 
 \caption{The five types of indecomposable non-split graphs, their paired degree sequences and distinguishing numbers.}
 \label{SplittableTable}
 \end{table}

  \begin{lemma}
 The paired degree sequences and distinguishing numbers of  the graphs in Table \ref{SplittableTable} are correct. 
 \end{lemma}

 \begin{proof}
 Again, we will just focus on verifying the distinguishing number of the graphs.  For a single-vertex graph, one color is clearly sufficient for a distinguishing labeling.  For $S(p,q)$, we computed its distinguishing number in  Section 2.1.  For $S_2(p_1, q_1, p_2, q_2, \cdots, p_m, q_m)$ with $m \ge 2$, no automorphism will map the vertices of $S(p_i, q_i)$ to that of $S(p_j, q_j)$ when $i \neq j$  because $p_i \neq p_j$.  Thus, it is enough for a labeling to be distinguishing for each $S(p_i, q_i)$.  That is, $D(S_2(p_1, q_1, p_2, q_2, \cdots, p_m, q_m)) = \max \{D(S(p_i, q_i)), i = 1, \hdots, m \}$. 
 
 For $S_3(p, q_1, q_2)$, every automorphism will fix vertex $e$ and, consequently, $S(p, q_1)$ and $S(p+1, q_2)$.  Thus, every distinguishing labeling will just have to be distinguishing for the two graphs so $D(S_3(p, q_1, q_2)) = \max \{D(S(p, q_1)), D(S(p+1, q_2)) \}$.  In $S_4(p, q)$, every automorphism will additionally fix vertex $f$ so $D(S_4(p, q)) = \max \{D(S(p, 2)), D(S(p+1, q)) \}$.
 \end{proof} 


  We emphasize that because of Proposition \ref{equaldistnums}, it is enough that we know the distinguishing numbers of the graphs in  Tables \ref{UnsplittableTable} and \ref{SplittableTable} to compute the distinguishing number of {\it any} indecomposable unigraph.  In the appendix, we present algorithms for computing  $D(mK_2)$ and $D(S(p,q))$ which we will then use to compute the distinguishing numbers of the other graphs in Tables \ref{UnsplittableTable} and \ref{SplittableTable}.  

\begin{lemma}
\label{distcomputelemma}
FindDistmK2($m$) computes $D(mK_2)$ with $m \ge 2$ in $O(\sqrt{m})$ time while FindDistS($p,q$) computes $D(S(p,q))$ with $p \ge 1, q \ge 2$ in $O(q)$ time.  
 \end{lemma}

The proof of Lemma \ref{distcomputelemma} can be found in the appendix. An immediate consequence of the lemma is  that  $D(mK_2)$ and $D(S(p,q))$ can be computed in time that is linear in the number of vertices of $mK_2$ and $S(p,q)$ respectively.

\section{Automorphisms and distinguishing labelings of $(G,A, B) \circ H$}

   We now establish several lemmas about the automorphisms and distinguishing labelings of $(G, A, B) \circ H$ so we can prove our first main result,  Theorem \ref{main1}.

   \begin{lemma}
   \label{autchar}
   Let $G^* = (G,A,B) \circ H$.  There is an automorphism $\pi$ of $G^*$ that maps $u \in V(G)$  to $v \in V(H)$ if and only if one of the following conditions hold: \smallskip
   
   (i) $u \in A$ and a swing vertex of $G$ while $v$ is a dominant vertex of $H$ or \smallskip
   
   (ii) $u \in B$ and a swing vertex of $G$ while $v$ is an isolated vertex of $H$.
  
   \end{lemma}
   
  
  \begin{proof} Assume $\pi \in Aut(G^*)$ such that $\pi(u) = v$ with $u \in V(G)$ and $v \in V(H)$.  Thus, $u \in A$ or $u \in B$.   When $u \in A$, $deg_{G^*}(u) \ge |A| - 1 + |V(H)|$ since $u$ is adjacent to all the vertices in $A$ (except itself) and in $H$ and possibly some vertices in $B$.  On the other hand, $deg_{G^*}(v) \leq |A| + |V(H)|  -1$ since $v$ is adjacent to all the vertices in $A$ and possibly all the vertices in $H$ (except itself) but none in $B$. In order for $\pi$ to map $u$ to $v$,  $deg_{G^*}(u) = deg_{G^*}(v) = |A| - 1 + |V(H)|$.   This means that $deg_G(u) = |A| -1$.  Since $u \in A$, by Proposition \ref{swingvertexprop}, $u$ is a swing vertex of $G$. On the other hand, in $H$, $v$ is adjacent to all the vertices in $H$ and is a dominant vertex.  Thus, (i) is true.

  
  Next, consider the case when $u \in B$.  Then $deg_{G^*}(u) \leq |A|$ because $u$ is not adjacent to any of the vertices in $B \cup V(H)$ while $deg_{G^*}(v) \ge |A|$ because $v$ is adjacent to all the vertices in $A$.  Again, if $\pi(u) = v$,  $deg_{G^*}(u) = deg_{G^*}(v) = |A|$.  Hence, $u \in B$ and $deg_G(u) = |A|$.  By Proposition \ref{swingvertexprop}, $u$ is a swing vertex of $G$. For $v$ though, $ deg_{G^*}(v) = |A|$ implies that it has no neighbors in $H$.  Condition (ii) is true.

   
 Let us now prove the converse.  Assume conditions (i) or (ii) holds.  Define the function $\eta: V(G^*) \rightarrow V(G^*)$ so that $\eta(u) = v$, $\eta(v) = u$ and $\eta(w) = w$ for $w \neq u, v$.  
  Pick any two vertices $x$ and $y$ in $V(G^*)$.  We shall show that $xy$ is an edge of $G^*$ if and only if $\eta(x) \eta(y)$ is also an edge of $G^*$.  This fact trivially holds when $x, y \not \in \{u, v\}$ or when $\{x, y\} = \{u, v\}$. 
  The interesting case is when $|\{x, y \} \cap \{u,v\} | = 1$. 
  
   Without loss of generality, let $x = u$ and $y \neq v$.  The case when $x = v$ and $y \neq u$ is argued similarly. When $u 
  \in A$ and a swing vertex of $G$,  $u$ is not adjacent to any vertices in $B$.  Thus,  in $G^*$,  $u$ is adjacent to the vertices in $A-\{u\} \cup V(H)$.  On the other hand, when $v$ is a dominant vertex in $H$,  $v$ is adjacent to $V(H) - \{v\} \cup A$ in $G^*$.   In other words, if condition (i) holds, $u$ and $v$ have the same neighbors in $V(G^*) - \{u, v\}$.   That is, $uy$ is an edge of $G^*$ if and only if $vy$ is an edge of $G^*$.   
  
  When $u \in B$ and a swing vertex of $G$, $u$ is adjacent to all the vertices in $A$, which are also its neighbors in $G^*$.  When $v$ is an isolated vertex in $H$, its neighbors in $G^*$ are exactly the vertices in $A$.   Thus, if condition (ii) holds, it is also the case that $uy$ is an edge of $G^*$ if and only if $vy$ is an edge of $G^*$.    We have now shown that when conditions (i) or (ii) hold, $\eta$ is an automorphism of $G^*$ that maps $u$ to $v$ and vice versa. 
  \end{proof}
  
  
  
  
  \begin{lemma}
  \label{autchar2}
  Let $G^* = (G,A,B) \circ H$.  Suppose one of the following conditions hold:
  \smallskip
  
(i) $(G,A,B)$ has no swing vertices in $A$ or $H$ has no dominant vertices;

(ii) $(G,A,B)$ has no swing vertices in $B$ or $H$ has no isolated vertices. 
\smallskip
  
\noindent Then a bijection $\pi: V(G^*) \rightarrow V(G^*)$ is an automorphism of $G^*$ if and only if $\pi$ when restricted to $V(G)$ is an automorphism of $(G,A,B)$ and $\pi$ when restricted to $V(H)$ is an automorphism of $H$. 

  
  \end{lemma}
 
 \begin{proof}  
 For ease of discussion, let $\pi_1$ and $\pi_2$ be equal to $\pi$ when restricted to $V(G)$ and $V(H)$ respectively.  Assume $\pi$ is an automorphism of $G^*$.  According to Lemma \ref{autchar}, if conditions (i) or (ii) holds,  $\pi$  fixes $V(G)$ and $V(H)$.   But in fact we can say something more about $\pi$.  Let $u \in A$. Since $u$ is adjacent to every vertex in $V(H)$, so must $\pi(u)$.  But vertices in $B$ are not adjacent to vertices in $H$.  So $\pi(u) \in A$. That is, $\pi$ fixes $A$ and, consequently, fixes $B$ too.   Now, the graphs induced by $V(G)$ and $V(H)$ in $G^*$ are $G$ and $H$ respectively.  Thus, $\pi_1$  and $\pi_2$ are automorphisms of $G$ and $H$. And since $\pi$ fixes $A$ and $B$, $\pi_1$ is also an automorphism of $(G,A,B)$. 
  
  Now, suppose $\pi_1$ is an automorphism of $(G,A,B)$ and $\pi_2$ is an automorphism of $H$.  To prove that $\pi$ is an automorphism of $G^*$,  we have to show that for any two vertices $x, y \in V(G^*)$, $xy$ is an edge of $G^*$ if and only if $\pi(x) \pi(y)$ is an edge of $G^*$.  But this statement is true when both $x,y \in V(G)$ or both $x, y \in V(H)$ since $\pi$ when restricted to $V(G)$ and $V(H)$ are just $\pi_1$ and $\pi_2$ respectively. 
  So assume $x \in V(G)$ and $y \in V(H)$; the case when $y \in V(G)$ and $x \in V(H)$ is addressed in the same way. Now, if $xy$ is an edge of $G^*$, then $x \in A$ because none of the vertices in $B$ are adjacent to the vertices of $H$.  But $\pi(x) = \pi_1(x) \in A$ since $\pi_1$ fixes $A$ and $\pi(y) = \pi_2(y) \in V(H)$.  Since every vertex in $A$ is adjacent to every vertex of $H$, $\pi(x) \pi(y)$ is an edge of $G^*$.   On the other hand, if $xy$ is not an edge of $G^*$, then $x \in B$.  But $\pi(x) = \pi_1(x)  \in B$ too.  It follows that  $\pi(x) \pi(y)$ is not an edge of $G^*$.  We have shown that $\pi$ is an automorphism of $G^*$.  \end{proof}


  
\begin{lemma}
\label{distlabelchar}
Let $G^* = (G,A,B) \circ H$. Suppose one of the following conditions hold:
  \smallskip
  
(i) $(G,A,B)$ has no swing vertices in $A$ or $H$ has no dominant vertices;

(ii) $(G,A,B)$ has no swing vertices in $B$ or $H$ has no isolated vertices. 
\smallskip
 
  Then $\phi$ is a distinguishing labeling of $G^*$ if and only if  $\phi$ is a distinguishing labeling of $(G,A,B)$ when restricted to $V(G)$ and 
$\phi$ is a distinguishing labeling of $H$ when restricted to $V(H)$. Consequently, $D(G^*) = \max\{D(G,A,B), D(H)\}$. 
\end{lemma}

\begin{proof}  
Assume $\phi$ is a distinguishing labeling of $G^*$ but $\phi$ does not distinguish $(G,A,B)$ when restricted to $V(G)$.  Then $(G,A,B)$ has a nontrivial automorphism $\eta$ so that $\phi(u) = \phi(\eta(u))$ for each $u \in V(G)$. Define $\overline{\eta}$ as the extension of $\eta$ to $G^*$ where $\overline{\eta}(u) = \eta(u)$ for each $u \in V(G)$ and  $\overline{\eta}(v) = v$ for each $v \in V(H)$.  By  Lemma \ref{autchar2},   $\overline{\eta}$ is a non-trivial automorphism of $G^*$.  Furthermore,  for every vertex $w$ of $G^*$, $\phi(w) = \phi(\overline{\eta}(w))$ so $\phi$ is not a distinguishing labeling of $G^*$,  a contradiction.  Thus, $\phi$ is a distinguishing labeling of $(G,A,B)$ when restricted to $V(G)$.  Using a similar argument, $\phi$ is a distinguishing labeling of $H$ when restricted to $V(H)$. 

 Next, consider the converse.  Suppose $\phi$ is a distinguishing labeling of $(G,A,B)$ when restricted to $V(G)$ and 
$\phi$ is a distinguishing labeling of $H$ when restricted to $V(H)$, but $\phi$ is not a distinguishing labeling of $G^*$.  Then there is some non-trivial automorphism $\pi$ of $G^*$ so that $\phi(u) = \phi(\pi(u))$ for each $u \in V(G^*)$.   According to Lemma \ref{autchar2}, $\pi$ when restricted to $V(G)$ is a non-trivial automorphism of $(G,A,B)$ or $\pi$ when restricted to $V(H)$ is a non-trivial automorphism of $H$.  Thus, it must mean that $\phi$ is not a distinguishing labeling of $(G,A,B)$ or $\phi$ is not a distinguishing labeling of $H$, which is a contradiction. So $\phi$ is also a distinguishing labeling of $G^*$.

We have established that every distinguishing labeling of $G^*$ is a combination of two distinguishing labelings: one for $(G,A,B)$ and another for $H$.  Thus, the number of colors used by such a labeling is at least  $\max\{D(G,A,B), D(H)\}$. But $V(G) \cap V(H) = \emptyset$,  so the colors used  for the vertices of $H$ can be reused for the vertices of $(G,A,B)$ and vice versa. It follows that  $\max\{D(G,A,B), D(H)\}$ colors are enough to distinguish both $(G,A,B)$ and $H$ and $D(G^*) = \max\{D(G,A,B), D(H)\}$. 
\end{proof}

We are now ready for our main result. \smallskip

 \begin{theorem}
 \label{main1}
 Let $G = (G_k, A_k, B_k) \circ (G_{k-1}, A_{k-1}, B_{k-1}) \circ \cdots \circ (G_1, A_1, B_1) \circ G_0$ be the compact canonical decomposition of $G$.  Then $\phi$ is a distinguishing labeling of $G$ if and only if $\phi$ is a distinguishing labeling of $G_i$ when restricted to $V(G_i)$ for $i = 0$ to $k$.  Furthermore, $D(G) = \max \{ D(G_i), i = 0, \hdots, k \}.$

   \end{theorem}

\begin{proof} We use induction on $k$ to prove the theorem.  When $k = 0$,  $G = G_0$ so the theorem is clearly true.  Assume it holds for all graphs whose compact canonical decompositions have $k$ components.  Let $G = (G_k, A_k, B_k) \circ H$ where $H = (G_{k-1}, A_{k-1}, B_{k-1}) \circ \cdots \circ (G_1, A_1, B_1) \circ G_0$.  There are three possibilities for $G_k$.  For each one, we shall argue that at least one of the conditions in Lemma \ref{distlabelchar} is satisfied. 
   
\smallskip

\noindent {\it Case 1: When $G_k$ is an indecomposable split graph with at least two vertices.} By  Proposition \ref{decomposable},   it does not have a swing vertex so conditions (i) and (ii) of Lemma \ref{distlabelchar} are satisfied.


\smallskip

\noindent {\it Case 2:  When $G_k$ is a complete graph.} That is, $G_k$ is itself a $K_1$-type trivial component or it was formed by maximally combining consecutive $K_1$-type trivial components of $G$.  
If $|V(H)| = 1$, then $H = G_0$ is a trivial component. Since $G_k$ is made up of $K_1$-type trivial components of $G$, $G_0$ should have been combined with them in the compact canonical decomposition of $G$.  Thus, $|V(H)| \ge 2$. 
 
 From Proposition  \ref{decomposable}, if $H$ has a dominant vertex, then $H$ is decomposable and the leftmost component is a trivial component of type $K_1$.  But the trivial component should have been combined with $G_k$ in the compact canonical decomposition of $G$.  Since this is a contradiction, $H$ has no dominant vertices and condition (i) of Lemma \ref{distlabelchar} is satisfied.

 


 
\smallskip

\noindent {\it Case 3: When $G_k$ is a graph on isolated vertices.} That is, $G_k$ is itself an $S_1$-type trivial component or is formed by maximally combining consecutive $S_1$-type trivial components of $G$. Again, if $|V(H)| = 1$, then $H$ should have been combined with the other trivial components that make up $G_k$ in the compact canonical decomposition of $G$.  So $|V(H)| \ge 2$. 

  From Proposition  \ref{decomposable},  if $H$ has an isolated vertex, then $H$ is decomposable and the leftmost component is a trivial component of type $S_1$.  Like the previous case, it implies that trivial components of the same type were not maximally combined in the compact canonical decomposition of $G$ -- a contradiction.  It follows that $H$ has no isolated vertices and condition (ii) of Lemma \ref{distlabelchar} is satisfied.

\smallskip

Since all three cases satisfy conditions (i) or (ii) of Lemma \ref{distlabelchar}, we conclude that $\phi$ is a distinguishing labeling of $G$ if and only if $\phi$ is a distinguishing labeling of $(G_k, A_k, B_k)$ and of $H$ when $\phi$ is restricted to $V(G_k)$ and $V(H)$ respectively.  Furthermore, $D(G) = \max\{D(G_k, A_k, B_k), D(H)\}$.   

But when $G_k$ is an indecomposable split graph with at least two vertices, it is a balanced split graph by  Proposition \ref{decomposable} so $Aut(G_k) = Aut(G_k, A_k, B_k)$ by Proposition \ref{equalautomorph}.  Therefore, $\phi$ is not only a distinguishing labeling of $(G_k, A_k, B_k)$, it is also a distinguishing labeling of $G_k$. Additionallly, $D(G_k, A_k, B_k) = D(G_k)$ by Proposition \ref{samedistprop2}. 

And when $G_k$ is a complete graph, $(A_k, B_k)$ is a $K$-max partition according to Definition \ref{compactdef}, so all the vertices of $G_k$ are in $A_k$ and $B_k$ is empty.  Similarly, when $G_k$ is a graph of isolated vertices, $(A_k, B_k)$ is an $S$-max partition according to Definition \ref{compactdef}.  For both cases, a distinguishing labeling of $(G_k, A_k, B_k)$ is also a distinguishing labeling of $G_k$ and $D(G_k, A_k, B_k) = |V(G_k)| = D(G_k)$.


Finally,  by the induction hypothesis, $\phi$ is a distinguishing labeling of $H$ if and only if $\phi$ is a distinguishing labeling of $G_i$ when restricted to $V(G_i)$ for $i = 0, \hdots, k-1$, and  $D(H) =    \max \{D(G_i), i = 0, \hdots, k-1\}$ so
 \begin{eqnarray*}
 D(G) & = & \max\{D(G_k, A_k, B_k), D(H)\} \\ 
       & = &   \max\{D(G_k),  \max \{ D(G_i),  i = 0, \hdots, k-1\} \}\\
       & = &  \max \{ D(G_i), i = 0, \hdots, k\}.
 \end{eqnarray*}
 The theorem holds by mathematical induction.
  \end{proof}

\begin{corollary}
 When $G$ is a threshold graph, $D(G)$ is equal to the size of the largest component in its compact canonical decomposition. 
\end{corollary}

\begin{proof} A threshold graph is a graph that can be constructed from an empty graph by repeatedly adding an isolated vertex or a dominating vertex.  Equivalently,  when $G$ is a threshold graph, every component in the canonical decomposition of $G$ is trivial.   Thus, every component in the compact canonical decomposition of $G$ is either a complete graph or a graph of isolated vertices.  For both cases, the distinguishing number of the component is equal to the size of the component.  The corollary follows from Theorem \ref{main1}.
\end{proof}

The graph in Example 3 is a threshold graph whose compact decomposition is $(1^2, \emptyset) \circ (\emptyset, 0^4)$.  It  has two components whose sizes are $2$ and $4$ respectively.  Thus, its distinguishing number is $4$. 

 \section{Computing the Distinguishing Numbers of Unigraphs}
 
  We now focus our attention on the algorithmic aspects of computing the distinguishing numbers of unigraphs.  We already know that a unigraph's compact canonical decomposition can be computed in linear time using the algorithm DECOMPOSE-C.  The next step is determine each component's distinguishing number and then output the largest among them.  Now there are three types of components that can arise:  a complete graph, a graph with isolated vertices or an indecomposable graph with at least two vertices.  Tyskevich has already identified all the possibilities for the last category and we have computed their distinguishing numbers in Section 2.   Thus, the only thing left to do is match each component $G_i$ to the appropriate graph. But because the graphs we are dealing with are all unigraphs, the matching problem is relatively simple -- it will just involve the comparison of degree sequences\footnote{All throughout this section, we will assume that the degree sequences are abbreviated. If a graph has degree sequence $(d_1^{r_1}, d_2^{r_2}, \hdots, d_s^{r_s})$, it will be stored as a list of ordered pairs: $((d_1, r_1), (d_2, r_2), \hdots, (d_s, r_s))$.  Thus, the length of the list equals the number of distinct vertex degrees the graph has.  A simple scan of the list reveals how many vertices have the highest degree, the second highest degree, etc. When the graph has a paired degree sequence $(d_1^{r_1}, \hdots, d_h^{r_h}; d_{h+1}^{r_{h+1}}, \hdots, d_s^{r_s})$, the information can be stored in two lists; the first list is for $(d_1^{r_1}, \hdots, d_h^{r_h})$ and the second list for $(d_{h+1}^{r_{h+1}}, \hdots, d_s^{r_s})$.}  because the degree sequence of a unigraph is unique up to isomorphism.    
  
  We leave it to the reader to verify the next proposition, keeping in mind that $\overline{(G,A,B)} = (\overline{G}, B, A)$, $(G,A,B)^I = (G^I, B, A)$ and $\overline{(G,A,B)^I} = (\overline{G^I}, A, B)$.

  
 
  \begin{proposition} 
  \label{degseqprop}
  The following are true:
  
 \noindent (i) For any graph $G$, the degree sequence of $\overline{G}$ can be obtained from the degree sequence of $G$ in $O(|V(G)|)$ time.  In particular, when $G$ has $s$ distinct vertex degrees, then so does $\overline{G}$.  
  \smallskip
    
\noindent  (ii)  If $G$ is split with a specific $KS$-partition $(A,B)$, then the paired degree sequences of $\overline{(G,A,B)}, (G,A,B)^I$ and $\overline{(G,A,B)^I }$ can each be obtained from the paired degree sequence of $(G,A,B)$ in $O(|V(G)|)$ time.    
 \end{proposition}

        In the next set of lemmas, we describe conditions that will allow us to quickly match an indecomposable unigraph  to one of the graphs in Table  \ref{UnsplittableTable} and \ref{SplittableTable} if it is indeed isomorphic to such a graph.

\begin{lemma}
\label{recognizeunsplit}
Let $G$ be an indecomposable non-split unigraph. 
The degree sequence of $G$ has the form 
\smallskip




\noindent $\bullet$  $(2^5)$ if and only if $G \cong C_5$.
\smallskip

\noindent $\bullet$ $(1^{r_1})$ if and only if $G \cong mK_2$ with $m = r_1/2$.
\smallskip

\noindent $\bullet$ $(d_1, 1^{r_2})$ with $r_2 > 1$ if and only if $G \cong U_2(m, \ell)$ with $m = (r_2 - d_1)/2$ and $\ell = d_1$.




\noindent $\bullet$ $(d_1, 2^{r_2})$ with $r_2 >1$ if and only if $G \cong U_3(m)$ with $m = (d_1 - 2)/2 $.  


\end{lemma}

\begin{proof}


From Table \ref{UnsplittableTable}, we can easily verify that if $G$ is isomorphic to one of the graphs mentioned in the table, the form of the degree sequence as stated in the lemma is correct.  So let us just prove the converse of each statement.

When the degree sequence of $G$ is $(2^5)$ or $(1^{r_1})$, it is obvious that $G$ is isomorphic to $C_5$ and $mK_1$ with $m = r_1/2$ respectively.   When the degree sequence of $G$ is of the form $(d_1^{r_1}, d_2^{r_2})$, there are four possibilities: the unigraph is isomorphic to $U_2(m, \ell)$ or its complement, $U_3(m)$ or its complement.  We note that the smallest degrees in $U_2(m, \ell)$ and $U_3(m)$ are $1$ and $2$ respectively and there are at least two vertices with the smallest degrees.  On the other hand, the complements of both $U_2(m, \ell)$ and $U_3(m)$ have a single node with the smallest degree.  Thus, when the degree sequence of $G$ is of the form $(d_1, 1^{r_2})$ with $r_2 > 1$, we can conclude that $G \cong U_2(m, \ell)$. The vertex with the largest degree is the center of $K_{1, \ell}$ and has degree $\ell$ so  $\ell = d_1$, while the degree-$1$ vertices that are not part of $K_{1, \ell}$ form $mK_2$ so   $m = (r_2-d_1)/2$. Similarly, when the degree sequence of $G$ is of the form $(d_1, 2^{r_2})$ with $r_2 > 1$, we conclude that $G \cong U_3(m)$.  The vertex with the highest degree is the center and has degree $2m + 2$.  So $m = (d_1 -2)/2$. \end{proof}

 \begin{lemma} 
 \label{recognizesplit}
 Let $G$ be an indecomposable split unigraph with at least two vertices. The paired degree sequence of $G$ has the form
 
\smallskip


\noindent $\bullet$  $(d_1^{r_1}; 1^{r_2})$ if and only if $G \cong S(p,q)$ with $p = r_2/r_1$ and $q = r_1$.

\smallskip

\noindent $\bullet$  $(d_1^{r_1}, \hdots, d_h^{r_h}; 1^{r_s})$ with $h \ge 2$ if and only if  $G \cong S_2(p_1, q_1, \hdots, p_h, q_h)$ with $p_i = d_i - N +1$  and $q_i = r_i$ for $i = 1$ to $h$,  and $N = \sum_{i=1}^h r_i$. 

\smallskip

\noindent $\bullet$   $(d_1^{r_1}; d_2, 1^{r_3})$ if and only if $G \cong S_3(p, q_1, q_2)$ with  $p = d_1 -  r_1$, $q_1 = d_2$ and $q_2 = r_1 - d_2$.

\smallskip

\noindent Finally, if  $G$ has a paired degree sequence of the form $(d_1, d_2^{r_2}; 2^{r_3})$ with $r_2 > 1$ and none of $\overline{G}, G^I, \overline{G^I}$ have degree sequences that conform to the ones described above, then $G \cong S_4(p,q)$ with $p = d_2 - r_2 - 1$ and $q = r_2 - 2$. 



\smallskip

\end{lemma}
 
 \begin{proof}
 Again, from Table \ref{SplittableTable}, we can easily verify that if $G$ is isomorphic to one of the graphs mentioned in the three bullet points, the form of the corresponding degree sequence is correct.  So we will just prove the converse of each statement.

 Let $(A,B)$ be the unique $KS$-partition of $G$.  When the paired degree sequence of $G$ is $(d_1^{r_1}; 1^{r_2})$, all  the vertices in $A$ have the same degree $d_1 > 1$ while all the vertices in $B$ have degree $1$.  If we consider only the edges between the $A$ and $B$, the  subgraph consists of stars whose centers are in $A$.  There are $r_1$ stars and each one has $r_2/ r_1$ leaves since the centers have to have the same degrees.  It follows that $G \cong S(p,q)$ with $p = r_2/r_1$ and $q = r_1$.

 
  When the paired degree sequence of $G$ is $(d_1^{r_1}, \hdots, d_h^{r_h}; 1^{r_s})$ with $h \ge 2$, the vertices in $A$ have $h$ different degrees while all the vertices in $B$ have degree $1$.  Again, the subgraph formed by the edges between $A$ and $B$ are stars whose centers are in $A$. But because there are $h$ different degrees in $A$, there has to be $h$ different kinds of stars as well.  Their sizes are clearly dependent on the $d_i$'s.  Let $N = \sum_{i=1}^h r_i$.  A vertex with degree $d_i$, $1 \leq i \leq h$, is the center of a star that has $d_i - (N-1)$ leaves  since it is also adjacent to all other vertices in $A$. Thus, $G \cong S_2(p_1, q_1, p_2, q_2, \hdots, p_h, q_h)$ where $p _i = d_i - N +1$, $q_i = r_i$ for $i = 1$ to $h$. 

 
 When the paired degree sequence of $G$ is  $(d_1^{r_1}; d_2, 1^{r_3})$, all the vertices in $A$ have the same degree $d_1$ while one vertex $v$ in $B$ has degree $d_2$ and the rest have degree $1$.   Let $A' \subseteq A$ contain the $d_2$ neighbors of $v$ in $A$. Notice that $A - A'$ cannot be empty; otherwise, $(G, A, B)$ is not a balanced split graph because $v$ will be a switch graph.   Let $a \in A'$ and $a'' \in A-A'$. Thus, the number of neighbors of $a $ in $B - \{v\}$ is one less than the number of neighbors of $a''$ in $B -\{v\}$.  Let the former number be $p$ so the latter number is $p+1$.  Since all the vertices in $B- \{v\}$ have degree $1$, it follows that the graph induced by $A' \cup (B-\{v\})$ is $S(p, d_2)$ while the graph induced by  $(A-A') \cup (B-\{v\})$ is $S(p+1, r_1 - d_2)$.  Since every vertex in $A$ has degree $d_1$ and it is adjacent to all other vertices in $A$ and $p+1$ vertices in $B$, $d_1 = r_1 - 1 + p + 1 = r_1 + p$.  Thus, $G \cong S_3(p, q_1, q_2)$  with  $p = d_1 -  r_1$, $q_1 = d_2$ and $q_2 = r_1 - d_2$.

 
 Lastly, if $G$'s paired degree sequence has the form $(d_1, d_2^{r_2}; 2^{r_3})$  with $r_2 > 1$ and none of its relatives are isomorphic to $S(p,q)$, $S(p_1, q_1, \hdots, p_h, q_h)$ or $S_3(p, q_1, q_2)$, then $G$ is isomorphic to $S_4(p,q)$, its complement, its inverse or the complement of its inverse.  We will now argue that the other three graphs are not valid options for $G$.
 
 Like $G$, the degree sequence of $S_4(p,q)$ has the form $(d_1, d_2^{r_2}; 2^{r_3})$ with $r_2 > 1$.  Two features stand out: there is only one vertex with the highest degree and all vertices in the $S$-part of the $KS$-partition have the same vertex degree $2$.  The $KS$-partitions of  the complement and the inverse of $S_4(p,q)$ is the reverse of the $KS$-partition of $S_4(p,q)$.  It is easy to check that their $S$-parts will have two different vertex degrees.  On the other hand, the complement of the inverse of $S_4(p,q)$ will have the same $KS$-partition as $S_4(p,q)$ but this time around there are $r_2 > 1$ vertices with the highest degree.  Thus, only one option is left for $G$ so $G \cong S_4(p,q)$.  There are $q+ 2$ vertices with the second largest degree so $q = r_2 - 2$.  The second largest degree is $p + q  + 3$ so $p = d_2 - q - 3 = d_2 - r_2  - 1$. \end{proof}

 \begin{remark}
 It would have been possible to create a fourth bullet point in Lemma \ref{recognizesplit} for $S_4(p,q)$ but the proof is significantly more involved.  We opted for the weaker final statement because it is enough to prove the correctness of the algorithm FindDistSplit presented below.
 \end{remark}
 


 \begin{figure}
\begin{algorithm}[H]
\caption{FindDistSplit($DegSeq$)}
\begin{algorithmic}[1]

\STATE Let $G$ be the split indecomposable unigraph with degree sequence $DegSeq$.
\STATE Let $Seq$ be an array so that $Seq[0] = DegSeq$ while $Seq[1]$, $Seq[2]$, $Seq[3]$ contain the degree sequences of $\overline{G}$, $G^I$ and $\overline{G^I}$ respectively. 
\STATE $distnum \leftarrow -1; i \leftarrow -1$
\IF{$G$ has two distinct degrees}
	\WHILE{$distnum = -1$  and $i < 3$}
		\STATE $i \leftarrow i +1$
			\IF{$S[i] = (d_1^{r_1}; 1^{r_2})$}
				\STATE $distnum \leftarrow$ FindDistS($r_2/r_1, r_1$)
				\RETURN $distnum$
			\ENDIF
	\ENDWHILE
\ELSIF{$G$ has three distinct degrees}
	\WHILE{$distnum = -1$ and $i < 3$ }
		\STATE $i \leftarrow i +1$
		\IF{$S[i] = (d_1^{r_1}, d_2^{r_2}; 1^{r_3})$}
			\STATE $distnum \leftarrow \max\{ \mbox{ FindDistS($d_1 - r_1 -r_2 + 1, r_1$),   FindDistS($d_2 - r_1 -r_2 + 1, r_2$)  } \}$
			\RETURN $distnum$
		\ENDIF
	\ENDWHILE
	\STATE $i = -1$
	\WHILE{$distnum = -1$ and $i < 3$ }
		\STATE $i \leftarrow i + 1$
		\IF{$S[i] = (d_1^{r_1}; d_2^{r_2}, 1^{r_3})$}
			\STATE $distnum \leftarrow \max\{ \mbox{ FindDistS($d_1 - r_1, d_2$),   FindDistS($d_1 - r_1 + 1, r_1 - d_2$)  } \}$
			\RETURN $distnum$
		\ENDIF
	\ENDWHILE
	\STATE $i = -1$
	\WHILE{$distnum = -1$ and $i < 3$ }
		\STATE $i \leftarrow i + 1$
		\IF{$S[i] = (d_1, d_2^{r_2}; 2^{r_3})$ and $r_2 > 1$}
			\STATE $distnum \leftarrow \max\{ \mbox{ FindDistS($d_2 - r_2 - 1, 2$),   FindDistS($d_2 - r_2, r_2 - 2$)  } \}$
			\RETURN $distnum$
		\ENDIF
	\ENDWHILE
\ELSE
	\STATE $i = -1$
	\WHILE{$distnum = -1$ and $i < 3$ }
		\STATE $i \leftarrow i +1$
		\STATE $(d_1^{r_1}, \hdots, d_h^{r_h}; d_{h+1}^{r_{h+1}}, \hdots,  d_s^{r_s}) \leftarrow S[i]$ 
		\IF{$s = h+1$ and $d_s = 1$}
			\STATE $N \leftarrow r_1 + r_2 + \hdots r_h$
			\FOR{$i = 1$ to $h$}
				\STATE $distnum \leftarrow \max\{ \mbox{FindDistS($d_i - N+1, r_i$)}, distnum\}$
			\ENDFOR
			\RETURN $distnum$
		\ENDIF
	\ENDWHILE
\ENDIF
\end{algorithmic}
\end{algorithm}
\label{finddistsplitalg}
\caption{A linear-time algorithm that determines the distinguishing number of an indecomposable split unigraph with at least two vertices.}
\end{figure}

 Next, we design the algorithm FindDistSplit in Figure \ref{finddistsplitalg}.  Given the paired degree sequence of an indecomposable split unigraph with size at least two, it computes the distinguishing number of the unigraph in linear time. 
 
 \begin{lemma}
 \label{runtime-splittable}
 Let $G$ be an indecomposable split unigraph with two or more vertices. Given the paired degree sequence of $G$, the algorithm FindDistSplit computes $D(G)$ in $O(|V(G)|)$ time. 
 \end{lemma}

\begin{proof}
From Theorem \ref{Tyshkevich-unsplittable}, we know that $G$ or one of its relatives has to be isomorphic to a graph in Table \ref{SplittableTable} because it is an indecomoposable split unigraph with at least two vertices.  We just have to figure out which one using Lemma \ref{recognizesplit}.   We then use the results in  Table \ref{SplittableTable} to compute $D(G)$.

 In FindDistSplit, the first step is to derive the degree sequences of all the relatives of $G$ from the degree sequence of $G$.   This step will take $O(|V(G)|)$ time according to Proposition \ref{degseqprop}.  These degree sequences are stored in the array $Seq[0 \hdots 3]$.  

Next, let $s$ be the number of distinct vertex degrees in $G$.  When $s = 2$,  some graph in $\{G, \overline{G}, G^I,   \overline{G^I}\}$ has to be isomorphic to $S(p,q)$. In lines 5 to 9, the algorithm checks which graph in the set has a degree sequence of the form $(d_1^{r_1}; 1^{r_2})$.  From Lemma \ref{recognizesplit}, $p = r_2/r_1$ and $q = r_1$ and 
FindDistS($r_2/r_1, r_1$) returns its distinguishing number.  Thus, after at most four iterations, the while loop will terminate and a number is returned.  This takes $O(|V(G)|)$.


When $s = 3$, it is harder to match $G$ or one of its relatives to a graph in Table \ref{SplittableTable} because there are three possibilities.  In lines 11 to 15, the algorithm goes through the graphs in $\{G, \overline{G}, G^I,   \overline{G^I}\}$ to see if one of them has 
a degree sequence of the form $(d_1^{r_1}, d_2^{r_2}; 1^{r_3})$. If so, then the graph is isomorphic to $S(p_1, q_1, p_2, q_2)$ with $p_i = d_i - r_1 - r_ 2 + 1$ and $q_i = r_i$ for $i = 1, 2$. 
Otherwise, in lines 16 to 21, the algorithm checks if one of the graphs 
has a degree sequence of the form $(d_1^{r_1}; d_2^{r_2}, 1^{r_3})$.  If yes, then it is isomorphic to $S_3(p, q_1, q_2)$ with $p = d_1 = r_1$, $q_1 = d_2$ and $q_2 = r_1 - d_2$.
When $distnum$ is still equal to $-1$ after line 21, it means neither $G$ nor any of its relatives were matched to $S(p,q)$, $S(p_1, q_1, \hdots, p_h, q_h)$ or $S_3(p,q_1, q_2)$.  Thus,  $G$ or one of its relatives has to be isomorphic to $S_4(p,q)$. In lines 22 to 27, the algorithm looks for a graph in $\{G, \overline{G}, G^I,   \overline{G^I}\}$ with a degree sequence of the form $(d_1, d_2^{r_2}; 2^{r_3})$ with $r_2 > 1$. Then $p = d_2 - r_2 -1$ and $q = r_2  - 2$. 


In the worst case, FindDistSplit may have to go through four iterations of the three while loops in lines 11-27 to identify the indecomposable split graph that matches $G$ or one of its relatives.  
When the right match is found, it computes the distinguishing number of the graph by calling   FindDistS($p,q$) twice and returning the maximum value.  Again, these set of steps take $O(|V(G)|)$ time.


Finally, when $s > 3$, $G$ or a relative  is isomorphic to $S(p_1, q_1, \hdots, p_h, q_h)$, $h \ge 3$.  Lines 20 to 37 goes through the four graphs to determine which one has  degree sequence of the form $(d_1^{r_1}, d_2^{r_2}, \hdots,  d_h^{r_h}; 1^{r_s})$.  Then $p_i = d_i - N+1$ and $q_i = r_i$ for $i = 1$ to $h$ and $N = \sum_{i=1}^h r_i$.  Once the match is made,  FindDistS($p_i,q_i$) is called for  $i = 1$ to $h$.  Checking the degree sequences take $O(|V(G)|)$ time and running the $h$ FindDistS calls also take $O(\sum_{i=1}^h q_i) = O(\sum_{i=1}^h r_i) = O(|V(G)|)$ time.   

We have shown that regardless of the value of $s$, FindDistSplit correctly computes the distinguishing number of indecomposable split graph $G$ in $O(|V(G)|)$ time. 
\end{proof}




We are now ready to present our algorithm FindDistUnigraph($G$), shown in Figure \ref{finddistunigraphalg},  that computes the distinguishing number of a unigraph $G$ in linear time.

 \begin{figure}
  \begin{algorithm}[H]
\caption{FindDistUnigraph($G$)}
\begin{algorithmic}[1]
  \STATE $T \leftarrow$ DECOMPOSE-C($G$)
  \STATE $distnum \leftarrow 1$
  \WHILE{$T.size > 1$ or ($T.size = 1$ and $T.top$ is a paired degree sequence)}
  	\STATE $DegSeq \leftarrow T.pop$ 
	\IF{$DegSeq = ((r_1-1)^{r_1}; \emptyset)$ or $DegSeq = (\emptyset; 0^{r_1})$}
		\STATE $distnum \leftarrow \max\{distnum, r_1\}$
	\ELSE
		\STATE $distnum \leftarrow \max \{ distnum, \mbox{FindDistSplit($DegSeq$)}\}$
	\ENDIF
  \ENDWHILE
  \IF{$T.size = 0$ or $T.top = (0)$}
  	\RETURN $distnum$
  \ELSE
        \STATE $DegSeq \leftarrow T.pop$
        \STATE $DeqSeq \leftarrow$ DetermineSplit($DegSeq$)
	\IF{$DegSeq$ is paired}
	         \RETURN $\max \{ distnum, \mbox{FindDistSplit($DegSeq$)}\}$
	\ELSE
  		\IF{$DegSeq$ has only one type of vertex degree}
			   \IF{ $DegSeq = (2^5)$}
			        \RETURN $\max \{ distnum, 3 \}$
			  \ELSIF{$DegSeq = (1^{r_1})$ or $DegSeq = ((r_1-2)^{r_1})$}
	 		 	 \RETURN $\max \{ distnum, \mbox{ FindDistmK2($r_1/2$)} \}$	
			   \ENDIF    
    		 \ELSE
			\STATE $(d_1^{r_1}, d_2^{r_2}) \leftarrow DegSeq$ 
			\IF{$r_1 = 1$}
				\IF{$d_2 = 1$}
					\STATE $c \leftarrow \max \{\mbox{FindDistmK2}((r_2 - d_1)/2), d_1\}$
				\ELSIF{$d_2 = 2$}	
					\STATE $c \leftarrow$ FindDistmK2($(d_1-2)/2$)
				\ENDIF		
			\ELSE
				\STATE $N = r_1 + r_2 $; $d'_1 = N-1 -d_2$; $d'_2 = N - 1 - d_1$; $r_1' = r_2$; $r_2' = r_1$
				\STATE $DegSeq \leftarrow ((d_1')^{r_1'}, (d_2')^{r_2'})$
				\IF{$d'_2 = 1$}
					\STATE $c \leftarrow \max \{\mbox{FindDistmK2}((r'_2 - d'_1)/2), d'_1\}$
				\ELSIF{$d'_2 = 2$}
					\STATE $c \leftarrow$ FindDistmK2($(d'_1-2)/2$)
				\ENDIF		
				\ENDIF			
			\ENDIF
		\RETURN $\max\{c, distnum\}$
		\ENDIF
      \ENDIF

  \end{algorithmic}
  \end{algorithm}
  \label{finddistunigraphalg}
  \caption{A algorithm that computes the distinguishing number of a unigraph in linear time.}
  \end{figure}

\begin{theorem}
 Let $G$ be a unigraph with $n$ vertices and $m$ edges.  FindDistUnigraph($G$) computes $D(G)$ in $O(n+m)$ time.  
  \end{theorem}

\begin{proof}
The algorithm begins by calling DECOMPOSE-C($G$) which returns a stack $T$ that contains the degree sequences of $G$'s compact canonical decomposition in order.  In particular, if $$G = (G_k, A_k, B_k) \circ (G_{k-1}, A_{k-1}, B_{k-1}) \circ \cdots \circ (G_1, A_1, B_1) \circ G_0$$ is the compact canonical decomposition of $G$, then the degree sequences of $G_k$ and $G_0$ are the top and bottom items of $T$.   According to Theorem \ref{decompose-c}, DECOMPOSE-C($G$) runs in  $O(n+m)$ time.



Initially, the algorithm sets $distnum$ to $1$.  While the top item $DegSeq$ in $T$ is a paired degree sequence, the algorithm removes $DegSeq$ and determines the distinguishing number of the associated graph $G_i$ in lines 3 to 8. We know that $G_i$ is either a complete graph, a graph of isolated vertices or an indecomposable split graph with two or more vertices. If it is the first two options, $D(G_i)$ is just $|V(G_i)|$; if the last option, the algorithm calls FindDistSplit to compute $D(G_i)$.  When $D(G_i) > distnum$, the algorithm updates $distnum$ to $D(G_i)$.  Thus, every iteration of the while loop takes $O(|V(G_i)|)$ time.


At the end of the while loop, $distnum = \max \{D(G_i), i = 0, \hdots, k\} = D(G)$ when the degree sequence of $G_0$ is paired so $D(G_0)$ was computed in the last iteration of the while loop. The algorithm returns $distnum$ in line 10.  But the other situation when we are sure that $distnum = D(G)$ is when $G_0$ is an isolated vertex because $D(G_0) = 1$ and will not cause any updates to $distnum$.  This case is also addressed in line 10.

 Otherwise,  $distnum = \max \{D(G_i), i = 1, \hdots, k\}$ and there is one graph left to consider -- $G_0$, which has at least two vertices. The fact that $G_0$'s degree sequence is not paired means that $G_0$ was not combined with any trivial components in the compact canonical decomposition of $G$.  Thus, $G_0$ may still be a split graph.  The algorithm checks this possibility in line 13; if yes, it again calls FindDistSplit to compute $D(G_0)$ and returns the larger of $distnum$ and $D(G_0)$ in line 15.   Otherwise $G_0$ is not a split graph. 

The rest of the algorithm is about identifying $G_0$ using Lemma \ref{recognizeunsplit} and computing its distinguishing number using Table \ref{UnsplittableTable}.  If $G_0$ has only one type of vertex degree, $G_0$ or its complement is  isomorphic to $C_5$ or $mK_2$.  It the former, $D(G_0) = 3$; if the latter, the algorithm calls FindDistmK2 to compute $D(G_0)$.  The algorithm returns the larger of $distnum$ and $D(G_0)$.  

On the other hand, if $G_0$ has two types of vertex degrees, the algorithm checks first if $G_0$ is isomorphic to $U_2(m, \ell)$ or $U_3(m)$ and computes its distinguishing number using FindDistmK2 in lines 24 to 28.  Otherwise, it set $DegSeq$ to the degree sequence of $\overline{G_0}$ and checks if $\overline{G_0}$ is the one that's isomorphic to $U_2(m, \ell)$ or  $U_3(m)$ and computes its distinguishing number in lines 30 to 35.  Again, the algorithm returns the larger of  $distnum$ and $D(G_0)$. 

It is straightforward to check that lines 23 to 36 also takes $O(|V(G_0)|)$ time since the algorithm makes only one call to FindmK2.  Thus, computing $D(G_0)$ takes $O(|V(G_0)|)$ time whether $G_0$ is a split graph or not.  Computing $D(G)$ then takes $\sum_{i=0}^h O(|V(G_i)| = O(|V(G)|)$ time.  \end{proof}

\noindent {\it Example 4:}    In \cite{BoCaPe11}, Borri et al. considered a $20$-vertex graph $G$ whose degree sequence is $(16^3, 12^4, 9^5, $ $5^2, 3, 2, 1^4)$.  Applying DECOMPOSE,  we get the following degree sequences stored in the stack $S$, from the bottom to the top:  $(4^3; 2, 1^4), (\emptyset; 0), (4^4; 2^2), (2^5)$.  In DECOMPOSE-C, none of the components are combined because there is only one trivial component.  
Finally, each of the degree sequences are processed to determine its type and distinguishing number.  

The degree sequence $(4^3; 2, 1^4)$ is paired and has the form $(d_1^{r_1}; d_2^{r_2}, 1^{r_3})$.  The graph is isomorphic to $S_3(1; 2, 1)$ so the algorithm computes $D(S(1,2))$ and $D(S(2,1))$ and outputs the maximum of the two values, which is $2$.  The graph with degree sequence $(\emptyset; 0)$ has distinguishing number $1$.   The degree sequence $(4^4; 2^2)$ is again paired but does not have the form $(d_1^{r_1}; 1^{r_2})$.  Its complement, however, has degree sequence $(3^2; 1^4)$.  So the complement of the graph is isomorphic to $S(2,2)$ and its distinguishing number is $2$.  Finally, $2^5$ is the degree sequence of $C_5$ whose distinguishing number is $3$.  Thus, $D(G) = 3$.




\section{Conclusion}

      In the introduction, we noted that unigraphs have a very simple isomorphism algorithm and we wanted to investigate if the algorithm can serve as the basis for designing an efficient algorithm for computing their distinguishing numbers.  Our work suggests that a more convoluted isomorphism algorithm for unigraphs is more useful:  given two unigraphs $G$ and $G'$, find their compact canonical decomposition using DECOMPOSE-C.  If they have the same number $r$ of components and the $i$th ones have identical degree sequences for $i = 0, \hdots, r-1$, then the two graphs are isomorphic; otherwise, they are not.   The algorithm still runs in time linear in the size of the two graphs.  Additionally, it reveals the structures of the unigraphs better because the components in the decomposition have limited type and, as we showed in the paper (Theorem \ref{main1}), the automorphisms of the unigraphs fix each component.  It can be thought of as the basis of our algorithm FindDistUnigraph. 
      
      We can also extend the result in Theorem \ref{main1} further. Recall that $D(G,c)$ is the number of inequivalent distinguishing labelings of $G$ that uses at most $c$ colors.

       
       
       
\begin{theorem}
Let  $G = (G_r, A_r, B_r) \circ (G_{r-1}, A_{r-1}, B_{r-1}) \circ \cdots \circ (G_1, A_1, B_1) \circ G_0$ be the compact canonical decomposition of $G$.  Then $D(G,c) = \prod_{i=0}^r D(G_i,c).$ 
\end{theorem}
      
      The theorem will be especially useful when multiple copies of $G$ have to be distinguished.  Finally, we pose the following question -- how might parameters related to distinguishing numbers like chromatic distinguishing numbers \cite{CoTr06}, list distinguishing numbers \cite{FeFlGe11}, etc. take advantage of Tyskevich's canonical decomposition?  Can an approach similar to ours lead to exact answers or good approximations?

        
      
\bibliography{dist-bib}

\bibliographystyle{abbrv}

 \newpage

\appendix

\section{Recognizing a split graph}

In Theorem \ref{splitthm}, Hammer and Simeone \cite{HaSi81} described how to determine if a graph is split using its degree sequence.  We use the theorem as a basis for the algorithm DetermineSplit($DegSeq$) in Figure \ref{DetermineSplitAlg}.  The input is the degree sequence of the graph.  If the graph is split, the algorithm transforms $DegSeq$ into a paired degree sequence that corresponds to a $KS$-partition of the graph; if not, $DegSeq$ stays the same.  The algorithm returns $DegSeq$. It is straightforward to check that DetermineSplit runs in time linear in the number of vertices of the graph.  
\medskip 

\noindent {\bf Corollary \ref{splitalgcor}.}  {\it Given the degree sequence of a graph with $n$ vertices, DetermineSplit determines if the graph is split in $O(n)$ time.}

  \begin{figure}[h]
\begin{algorithm}[H]
\caption{DetermineSplit($DegSeq$)}
\begin{algorithmic}[1]
 \STATE $(d_1, d_2, \hdots, d_n) \leftarrow DegSeq$
 \STATE $i \leftarrow 1$
 \WHILE{ $i \leq n-1$ and $d_{i+1} \ge i$}
     \STATE $i \leftarrow i+1$
 \ENDWHILE
 \STATE $h \leftarrow i$
 \STATE $leftsum \leftarrow 0$; $rightsum \leftarrow 0$
 \FOR{$i = 1$ to $h$}
    \STATE $leftsum \leftarrow leftsum + d_i$
  \ENDFOR
  \FOR{$i = h+1$ to $n$}
    \STATE $rightsum \leftarrow rightsum + d_i$
   \ENDFOR
   \IF{$leftsum = h(h-1) + rightsum$}
      \STATE $DegSeq \leftarrow (d_1,  \hdots, d_h; d_{h+1}, \hdots, d_n)$
    \ENDIF
  \RETURN($DegSeq$)
\end{algorithmic}
\end{algorithm}
\caption{Given the degree sequence of a graph, the algorithm outputs a paired degree sequence if the graph is split; otherwise, it just outputs the original degree sequence.}
\label{DetermineSplitAlg}
\end{figure}

\section{Algorithms for the Canonical and Compact Canonical Decomposition of a Graph}

        In \cite{Ty00}, Tyshkevich alluded to a linear-time algorithm that computes the canonical decomposition of a graph.  She laid the groundwork but did not actually present the algorithm so we do it here. Our initial version of the algorithm is less efficient but easier to follow.   We then describe a modification so it runs in linear time.  

		
\begin{theorem} (\cite{Ty00})
\label{goodpair}
An $n$-vertex graph $G$ with a degree sequence $(d_1, d_2, \hdots, d_n)$ is decomposable if and only if there exists non-negative integers $p$ and $q$ such that $$ 0 < p + q < n \; \; \;{(\star)} $$ and $$ \sum_{i = 1}^p d_i = p(n-q-1) + \sum_{i = n-q+1}^{n} d_i.  \; \; \; \;{(\star \star)}  $$
Call such a pair $(p,q)$ good.  For every good pair $(p,q)$, the decomposition is $G = (G', A', B') \circ H$ where 
$$(d_1, \hdots, d_p), (d_{p+1}, \hdots, d_{n-q}) \mbox{ and } (d_{n-q+1}, \hdots, d_n)  $$ 
are the degree sequences for the vertices in $A'$, $H$ and $B'$ respectively.  Moreover, every decomposition of $G$ is associated with some good pair $(p,q)$. 
\end{theorem}

\begin{remark}
\label{tyskevichremark}
We emphasize that when $(p,q)$ is a good pair of $G$, $|A'| = p$ and $|B'| = q$ so $G'$ has $p+q$ vertices while $H$ has $n-p-q$ vertices.  Moreover, the vertices of $A'$ have the $p$ largest degrees in $G$ while the vertices of $B'$ have the smallest $q$ degrees in $G$.  The paired degree sequence of $(G', A', B')$ then is 
$(d_1 - \alpha, d_2 - \alpha, \hdots, d_p - \alpha;  d_{n-q+1}, \hdots, d_n) $
with $\alpha = (n-p-q)$ because every vertex in $A'$ is adjacent to all the vertices of $H$ while none of the vertices in $B'$ have neighbors in $H$.  For the same reason, the degree sequence of $H$ is $(d_{p+1} - \beta, \hdots, d_{n-q} - \beta)$
with $\beta = p$.   We shall refer to $\alpha$ and $\beta$ as adjustment values.
\end{remark}

\medskip
According to Theorem \ref{goodpair}, if the canonical decomposition of  $G$ is $$G = (G_r, A_r, B_r) \circ (G_{r-1}, A_{r-1}, B_{r-1}) \circ \cdots \circ (G_1, A_1, B_1) \circ G_0,$$ then $G$ has $r$ good pairs.  Of interest to us is the good pair that separates $(G_r, A_r, B_r)$ from the rest of $G$. Notice that $(G_r, A_r, B_r)$ must be an indecomposable split graph.  Tyskevich described how to find this specific good pair.

\begin{corollary} (\cite{Ty00})
\label{goodpaircor}
Let the canonical decomposition of $G$ be 
$$G = (G_r, A_r, B_r) \circ (G_{r-1}, A_{r-1}, B_{r-1}) \circ \cdots \circ (G_1, A_1, B_1) \circ G_0.$$
The good pair  that corresponds to the decomposition $G = (G_r, A_r, B_r) \circ H$ where $H = (G_{r-1}, A_{r-1}, B_{r-1}) \circ \cdots \circ (G_1, A_1, B_1) \circ G_0$  is the lexicographically least among all good pairs of $G$.  Specifically, it is either $$ (0,1) \mbox{ or has the form } (p, |\{i: d_i < p\}|) \mbox{ with $p \ge 1$.} \;\;\; {(\star \star \star)} $$

\end{corollary}

\begin{remark}
\label{goodpairremark}
Since we are interested in finding the lexicographically least good pair of $G$,  we start our search with $(0,1)$.  According to ($\star \star)$,  $(0,1)$ is a good pair if and only if $d_n = 0$; i.e., when $G$ has an isolated vertex. If $(0,1)$ is not a good pair, then all the vertices of $G$ have degrees at least $1$.  Next, we consider the pair $(1,0)$ because when  $p = 1$,  $q = |\{i :  d_i = 0\}| = 0$ by  ($\star \star \star$).  Again, according to ($\star \star)$, $(1,0)$ is a good pair if and only if $d_1 = n-1$; i.e., when $G$ has a dominating vertex. 


\end{remark}


Our main algorithm for Tyskevich's decomposition is DECOMPOSE($G$), shown in Figure \ref{hello0}. It returns a stack $S$ that contains the (paired) degree sequences of the indecomposable components of $G$ with $G_0$ at the top.  As indecomposable components are peeled off from $G$, a smaller graph $H$ is left whose degree sequence in $G$ is $(d_{i+1}, d_{i+2}, \hdots, d_j)$ for some $i$ and $j$ with $i < j$.  But the true degree sequence of $H$,  $(d'_1, d'_2, \hdots, d'_{j-i})$, is not necessarily the same\footnote{The only exception is when the indecomposable components that were removed from $G$ are all isolated vertices.} as noted in Remark \ref{tyskevichremark}.   The algorithm uses $\beta$ to keep track of the adjustment value; i.e. $d'_t = d_{i+ t} - \beta$, for $t = 1, 2, \hdots, j-i$.

The next step is to determine if $H$ has a good pair using the subroutine FindGoodPair, shown in Figure \ref{hello}.  In particular, FindGoodPair uses conditions $(\star), (\star \star)$ and $(\star \star \star)$ to identify such a pair. 
 If $H$ has no good pairs then the graph is indecomposable and the algorithm has completed the canonical decomposition of $G$.  The degree sequence of $H$ is added to $S$,  and $S$ is returned.  But if $H$ has a good pair $(p,q)$ then the paired degree sequence of the indecomposable component is added to $S$.  Once the indecomposable component is removed from $H$, the remaining graph's degree sequence in $G$ is $(d_{i+p+1}, \hdots, d_{j-q})$ so $i$ is incremented by $p$ while $j$ is decremented by $q$.  The adjustment value $\beta$ is also increased by $p$.

 \begin{figure}[h]
\begin{algorithm}[H]
\caption{DECOMPOSE($G$)}
\begin{algorithmic}[1]
\STATE Let $D[1 \cdots n]$ contain the degree sequence of $G$ with $D[1] \ge D[2] \ge \cdots \ge D[n]$. 
\STATE $i \leftarrow 0$; $j \leftarrow n$
\STATE $\beta \leftarrow 0$
\STATE Let $S$ be an empty stack.
\WHILE{$i < j$}
	\IF{$j-i = 1$}
		\STATE $S.push((0))$
		\RETURN $S$
	\ENDIF
	\STATE $(p, q) \leftarrow$ FindGoodPair($D, i, j,  \beta$)
	\IF{$(p,q) = (-1, -1)$}
		\STATE $S.push((D[i+1] - \beta, \hdots, D[j] - \beta))$
		\RETURN $S$
	\ELSE	
		\IF{$(p,q) = (0,1)$}
			\STATE $S.push((\emptyset ; 0))$
		\ELSIF{$(p,q) = (1,0)$}
			\STATE $S.push((0; \emptyset))$
		\ELSE
			\STATE $\alpha \leftarrow j - i - p - q$
			\STATE $S.push((D[i+1] - \beta - \alpha, \hdots, D[i+p] - \beta - \alpha; D[j -q +1] - \beta, \hdots, D[j] - \beta))$
		\ENDIF
		\STATE $i \leftarrow i+p$; $j \leftarrow j-q$
		\STATE $\beta \leftarrow \beta + p$
         \ENDIF
\ENDWHILE
\end{algorithmic}
\end{algorithm}
\caption{An algorithm that outputs the canonical decomposition of a graph $G$.}
\label{hello0}
\end{figure}

 \begin{figure}[h]	
\begin{algorithm}[H]
\caption{FindGoodPair($D, i, j, \beta$)}
\begin{algorithmic}[1]
	\STATE \COMMENT{We shall use $m$  for the number of vertices in the graph and array $D'$ to contain its degree sequence.}
	\STATE $m \leftarrow j - i$  
	\STATE Create array $D'[1 \cdots m]$.
	\FOR{$t = 1$ to $m$}
		\STATE $D'[t] \leftarrow D[i+t] - \beta$ 
	\ENDFOR
	\IF{$D'[m] =  0$}
		\RETURN $(0,1)$
	\ELSIF{$D'[1] = m-1$}
		\RETURN $(1,0)$
	\ELSE
		\STATE $p \leftarrow 1$; $q = 0$
		\STATE $frontsum \leftarrow D'[1]$; $backsum \leftarrow 0$
		\WHILE{$p + q < m $ and $frontsum \neq p(m-q-1) + backsum$}
			\STATE $p \leftarrow p+1$
			\STATE $frontsum \leftarrow frontsum + D'[p]$
			\WHILE{$p + q < m $ and $D'[m-q] < p$}
				\STATE $q \leftarrow q+1$
				\STATE $backsum \leftarrow backsum + D'[m-q + 1]$
			\ENDWHILE
		\ENDWHILE
		\IF{$p+q < m$ and $frontsum = p(m-q-1) + backsum$}
			\RETURN $(p,q)$
		\ELSE
			\RETURN $(-1, -1)$
		\ENDIF
	\ENDIF
\end{algorithmic}
\end{algorithm}
\caption{An algorithm that finds the lexicographically least good pair in a graph whose degree sequence is $(D[i+1] - \beta, D[i+2] - \beta, \hdots, D[j] - \beta)$. If lines 3 and 4 are ignored and, for any integer $t$, $D'[t]$ is replaced by $D[i+t] - \beta$ throughout the pseudocode, the running time of the algorithm is $O(p + q)$ when it returns $(p,q)$  and $O(m) = O(j - i)$ when it returns $(-1, -1)$.}
\label{hello}
\end{figure}

 \begin{theorem}
  The algorithm DECOMPOSE($G$) correctly outputs the canonical decomposition of $G$. 
 \end{theorem}
 
 \begin{proof}
     Let us begin by addressing the correctness of FindGoodPair.  Let $H$ be the graph whose degree sequence is stored in $D'[1 \cdots m]$.  Assume $H$ has a good pair and the lexicographically least one is $(p_0, q_0)$.  One possibility for $(p_0, q_0)$ is $(0,1)$, and another is $(1,0)$.  The algorithm uses Remark \ref{goodpairremark} to check for these possibilities in lines 6 and 8 respectively.  Otherwise, the algorithm sets $p = 1$, $q = 0$, $frontsum = D'[1]$ and $backsum = 0$.  We already know that $(1,0)$ is not a good pair so the algorithm enters the while loop in lines 13 to 18.

     
      At each iteration of the while loop, the algorithm increments $p$ by $1$ and updates $frontsum$ so it is equal to $\sum_{t = 1}^p D'[t]$.  It then searches the right $q$-partner for $p$ using ($\star \star \star$).  Before entering the second while loop in lines 16 to 18,  the algorithm knows that the smallest $q$ degrees $D'[m], D'[m-1], \hdots, D'[m-q+1]$ are all less than $p-1$ and therefore less than $p$ too.  It goes through the second while loop to  determine if there are any additional degrees less than $p$ starting with $D'[m-q]$. A check is made to ensure that ($\star$) still holds.   If so, the algorithm increments $q$  and updates $backsum$ so it is equal to $\sum_{t = m-q +1}^m D'[t]$.

     Since $(p_0, q_0)$ is the lexicographically least good pair of $H$,  the algorithm keeps updating $p$ until $p = p_0$.  It then updates $q$ to $q_0$ and exits the outer while loop.  It does one last check that $(p_0, q_0)$ is a good pair  in line 19 and returns the pair.  
     
     When $H$ has no good pairs, the algorithm will keep incrementing $p$ and find its appropriate $q$-partner.  Eventually, $p+q \ge m$, violating  ($\star$), so the algorithm returns $(-1,-1)$ to indicate $H$ has no good pairs. 
     
     Now that we have established that FindGoodPair finds the lexicographically least good pair of $H$ if one exists, let us prove the correctness of DECOMPOSE. It works by peeling off the indecomposable components of $G$ and storing their (paired) degree sequences in the stack $S$.  At the beginning of the while loop in line 5, the remaining graph $H$ has degree sequence $(D[i+1] - \beta, D[i+2] - \beta, \hdots, D[j]-\beta)$.  When $j-i = 1$, $H$ has only one vertex and is therefore indecomposable.  The algorithm adds its degree sequence $(0)$ to $S$ and returns $S$.  The algorithm does the same thing when $H$ has no good pairs because it means $H$ is indecomposable.  However, when $H$ is decomposable,  the indecomposable component associated with its lexicographically least good pair $(p,q)$ has three possibilities: it is either an isolated or a dominating vertex in $H$ or neither.  The paired degree sequence of the third case is described in Remark \ref{tyskevichremark}.  The algorithm adds the appropriate paired degree sequence to $S$ and updates $i, j$ and $\beta$ accordingly.  \end{proof}

  Let $G$ be a graph with $n$ vertices and $m$ edges, and its canonical decomposition is 
  $$G = (G_r, A_r, B_r) \circ (G_{r-1}, A_{r-1}, B_{r-1}) \circ \cdots \circ (G_1, A_1, B_1) \circ G_0.$$
Assume the adjacency list of $G$ is given.  We now analyze the running time of DECOMPOSE($G$).  
 The first step of the algorithm  is to compute the degree sequence of $G$.  Obtaining the vertex degrees of $G$ takes $O(n+m)$ time.  Sorting them from largest to smallest using bucket sort takes $O(n)$ time.  Next, the algorithm enters the while loop in lines 5 to 22.  At every iteration, an indecomposable component $G_k$ of $G$ is peeled off and the degree sequence of $G_k$ is added to $S$.  Ignoring the call to FindGoodPair for the moment,  it is easy to verify that all the other steps take $O(|V(G_k)|)$ time.  Thus, if FindGoodPair also takes $O(|V(G_k)|)$ time, then the total running time of the while loop is $O(\sum_{k=0}^r |V(G_k)|) = O(n)$. Consequently, DECOMPOSE($G$) runs in $O(n+m)$ time. 
 
 \smallskip
  
The goal of  FindGoodPair($D,i,j,\beta$)  is to find a good pair for the graph $H$ whose degree sequence is $(D[i+1] - \beta, D[i+2] - \beta, \hdots, D[j] - \beta)$.  If it returns $(-1, -1)$, $H$ has no good pair and the indecomposable component that is added to $S$ is actually $H$ itself, which has $m = j-i$ vertices.  However, if it returns $(p,q) \neq (-1, -1)$, $H$ has a good pair and the indecomposable component that is added to $S$ has $p + q$ vertices, where $p+q < m$ in this case. 
 
 So consider FindGoodPair($D,i,j,\beta$).   It recalibrates the degree sequences of all the vertices and stores them in array $D'$ in lines 4 and 5. This step takes $O(m)$ time.  It then checks if $(0,1)$ or $(1,0)$ are good pairs in $O(1)$ time.  Otherwise, the algorithm searches for a good pair in lines 13 to 18, incrementing $p$ and $q$ as needed.
  Consequently, if FindGoodPair returns $(p,q) \neq (-1, -1)$, the running time of the while loop is $O(p+q)$.  However, if it returns $(-1, -1)$, the running time of the while loop is $O(m)$ because $p$ and $q$ are incremented until $p+q = m$.  
  
  Our running time analysis shows that regardless of what FindGoodPair returns, it will run in $O(m)$ time because of lines 4 and 5.  But in fact, we do not need these steps.  We recalibrated the degree sequences of the vertices so that indices of the degrees match those in Theorem \ref{goodpair} and Corollary \ref{goodpaircor}, making it easier to verify the correctness of the algorithm.  Instead, for any integer $t$, every time $D'[t]$ appears in the pseudocode, simply replace it with $D[i+t] - \beta$.  With this simple modification, FindGoodPair runs in $O(|V(G_k)|)$ time when $G_k$ is the indecomposable component that is being peeled off from $G$.

 \medskip
 
 \noindent  {\bf Theorem \ref{decomposealgthm}.}  {\it  Let $G$ be a graph with $n$ vertices and $m$ edges.  
 Suppose the canonical decomposition of $G$ is  $G = (G_r, A_r, B_r) \circ (G_{r-1}, A_{r-1}, B_{r-1}) \circ \cdots \circ (G_1, A_1, B_1) \circ G_0.$
The algorithm DECOMPOSE($G$) returns a stack $S$ that contains the (paired) degree sequences of the $G_i$'s in order  in $O(n+m)$ time.}

\medskip
 
 


Next, we present DECOMPOSE-C($G$) in Figure \ref{decompose-calg} which takes the output of DECOMPOSE($G$) and maximally combines the trivial components of the same type to produce the compact canonical decomposition of $G$.  Initially, the degree sequence of $G_0$ is compared with that of $G_1$ to see if the two components can be combined.  From there, every trivial component of $G$ is examined to see if it can be combined with the previously processed component.


\medskip

\noindent {\bf Theorem \ref{decompose-c}} {\it Let $G$ be a graph with $n$ vertices and $m$ edges.  
Suppose the compact canonical decomposition of $G$ is  $G = (G_k, A_k, B_k) \circ (G_{k-1}, A_{k-1}, B_{k-1}) \circ \cdots \circ (G_1, A_1, B_1) \circ G_0.$
The algorithm DECOMPOSE-C($G$) returns a stack $T$ that contains the abbreviated (paired) degree sequences of the $G_i$'s in order in $O(n+m)$ time. }


\begin{figure}
\begin{algorithm}[H]
\caption{DECOMPOSE-C($G$)}
\begin{algorithmic}[1]
\STATE $S \leftarrow$ DECOMPOSE($G$)
\STATE Go through all the degree sequences stored in $G$ and abbreviate them.
\IF{$S.size = 1$}
	\RETURN $S$
\ELSE
	\STATE Let $T$ be an empty stack. 
	\STATE $TopDegSeq \leftarrow S.pop$
	\IF{$TopDegSeq = (0)$ and $S.top = (0; \emptyset)$}
		\STATE $S.pop$
		\STATE $T.push((1^2; \emptyset))$
	\ELSIF{$TopDegSeq = (0)$ and $S.top = (\emptyset; 0)$}
		\STATE $S.pop$
		\STATE $T.push((\emptyset; 0^2))$
	\ELSE
		\STATE $T.push(TopDegSeq))$
	\ENDIF
	\WHILE{$S$ is not empty}
		\STATE $DegSeq \leftarrow S.pop$
		\IF{$DegSeq = (0; \emptyset)$ and $T.top = ((m-2)^{m-1}; \emptyset)$}
			\STATE $T.pop$
			\STATE $T.push(((m-1)^{m}; \emptyset))$
		\ELSIF{$DegSeq = (\emptyset; 0)$ and $T.top = (\emptyset; 0^{m-1})$}
			\STATE $T.pop$
			\STATE $T.push((\emptyset; 0^m))$
		\ELSE
			\STATE $T.push(DegSeq))$
		\ENDIF	
	\ENDWHILE
	\RETURN $T$
\ENDIF
\end{algorithmic}
\end{algorithm}
\caption{An algorithm that computes the compact canonical decomposition of a graph.}
\label{decompose-calg}
\end{figure}

  \section{Computing $D(mK_2)$ and $D(S(p,q))$}
   
According to Tables \ref{UnsplittableTable} and \ref{SplittableTable},  $D(mK_2) = \min \{c: \binom{c}{2} \ge m\}$ while $D(S(p,q)) = \min \{c: c\binom{c}{p} \ge q \}$.  FindDistmK2($m$) and FindDistS($p,q$) in Figure \ref{finddist-splitalg} use these facts to compute the distinguishing numbers of $mK_2$ and $S(p,q)$.  Their running times are linear in the number of vertices of the graph.

 \begin{figure}[h]
\label{finddistalg}
\begin{algorithm}[H]
\caption{FindDistmK2($m$)}
\begin{algorithmic}[1]
\STATE $curr \leftarrow 2$
\STATE $val \leftarrow 1$
\WHILE{$val < m$}
    \STATE $curr \leftarrow curr+1$
  \STATE $val \leftarrow val + curr -1$
\ENDWHILE
\RETURN $curr$
\end{algorithmic}
\end{algorithm}

\begin{algorithm}[H]
\caption{FindDistS($p,q$)}
\begin{algorithmic}[1]
\STATE $curr \leftarrow p$
\STATE $val \leftarrow p$
\WHILE{$val < q$}
    \STATE $curr \leftarrow curr+1$
  \STATE $val \leftarrow val\times \frac{(curr)^2}{(curr-1)(curr -p)}$
\ENDWHILE
\RETURN $curr$
\end{algorithmic}
\end{algorithm}
\caption{An $O(\sqrt{m})$-time algorithm that outputs $D(mK_2)$ and an $O(q)$-time algorithm that outputs $D(S(p,q))$.}
\label{finddist-splitalg}

\end{figure}

   \medskip

\noindent {\bf Lemma \ref{distcomputelemma}.} {\it
FindDistmK2($m$) computes $D(mK_2)$ with $m \ge 2$ in $O(\sqrt{m})$ time while FindDistS($p,q$) computes $D(S(p,q))$ with $p \ge 1, q \ge 2$ in $O(q)$ time.}  
\smallskip

\begin{proof}
 We already know that $D(mK_2) = \min \{c: \binom{c}{2} \ge m\}$.  Analytically, all we have to do is find the positive root of the quadratic equation $c(c-1) - 2m = 0$ and take its ceiling.  Thus, $D(mK_2) =  \lceil (1 + \sqrt{1+ 8m})/2\rceil$.  But we would like to avoid taking square roots because such operations do not run in constant time.

 Instead, in FindDistmK2($m$),  we make use of the observation that $\binom{c+1}{2} - \binom{c}{2} = c$.
 Thus, to compute $D(mK_2)$, let  $curr$ equal  the current $c$-value and $val = \binom{curr}{2}$.  Initially, set $curr$ to $2$ and $val$ to $1$.  While $val < m$, increase $curr$ by $1$ and update $val$ to $val + curr-1$. Once $val \ge m$, we know $D(mK_2) = curr$.   Each iteration takes $O(1)$ time and the number of iterations is $D(mK_2) -2$.   But $D(mK_2) \leq \sqrt{2m} \rceil + 1$ because  $(\sqrt{2m} +1)(\sqrt{2m})/2 = m + \sqrt{2m}/2 > m$. Thus, $D(mK_2)$ can be computed in $O(\sqrt{m})$ time.  
 


Next, consider $D(S(p,q))$, which is equal to $\min \{c: c\binom{c}{p} \ge q \}$.   Note that $D(S(p,q)) \ge p$ since $\binom{c}{p}$ has to be at least $1$.  Furthermore, when $p \ge q$,  $p$ colors are enough to distinguish $S(p,q)$ because $p \binom{p}{p} = p \ge q$.   On the other hand, when $q > p$, $q$ colors are enough because $q \binom{q}{p} >  q$. So $p \leq D(G_0) \leq \max \{p,q\}$.  
 
In the algorithm FindDistS($p,q$), we use the same idea for computing $D(mK_2)$  to compute $D(S(p,q))$.  The variables $curr$ and $val$ denote the current $c$-value and $c \binom{c}{p}$ respectively. When it's time to update $val$, instead of computing $(c+1) \binom{c+1}{p}$ from scratch, we make use of the fact that 
 $$(c+1) \binom{c+1}{p} - c\binom{c}{p} = \binom{c}{p} \left[ \frac{(c+1)^2}{c+1-p} - c\right]. $$

 
 Thus, when  $curr = c+1$ and $val = c \binom{c}{p}$, 
 \begin{eqnarray*}
 (c+1) \binom{c+1}{p} & = & val + \frac{val}{c} \left[ \frac{(c+1)^2}{c+1-p} - c\right] \\
   & = & val \left[1 + \frac{(c+1)^2}{c(c+1-p)} - 1 \right] \\
   & = &  val \left[\frac{(c+1)^2}{c(c+1-p)} \right] \\
   & = & val \left[\frac{curr^2}{(curr-1)(curr-p)} \right].
 \end{eqnarray*}
 Hence, every iteration of the while loop just takes $O(1)$ time.  The number of iterations of the while loop  is $0$ when $p \ge q$ and at most $q-p$ when $q > p$.  Thus, FindDistS($p,q$) computes $D(S(p,q))$ in $O(\max\{1, q-p \}) = O(q)$ time.  \end{proof}

\end{document}